\documentclass[final, 3p, times]{elsarticle}
\usepackage{amsmath, amssymb}
\usepackage{graphicx}
\usepackage[usenames]{color}
\usepackage{amsfonts}
\usepackage{dcolumn}
\usepackage{bm}%
\usepackage{pifont}
\usepackage{marvosym}
\usepackage[colorlinks=true,linkcolor=blue]{hyperref}%
\usepackage[mathlines]{lineno}
\usepackage{textgreek}
\usepackage{wasysym}
\usepackage{xspace}

\usepackage{etoolbox}

\newcommand{\be}{\begin{equation}}
\newcommand{\en}{\end{equation}}

\newcommand{\lshad}{[\![}
\newcommand{\rshad}{]\!]}
\newcommand{\rd}{\mathrm{d}}

\newcommand{\ram}{\text{\Zodiac{1}}}
\newcommand{\bull}{\text{\Zodiac{2}}}

\DeclareMathOperator{\sgn}{sgn}
\DeclareMathOperator{\avg}{\text{avg}}

\AtBeginDocument{\mathcode`v=\varv}

\biboptions{sort&compress}

\begin{document}

\title{Acoustic singular surfaces in an exponential class of inhomogeneous gases:  A new numerical approach based on Krylov subspace spectral methodologies\tnoteref{T1}}

\tnotetext[T1]{Distribution Statement A: Approved for public release. Distribution is unlimited.}
%
%

\author[USM]{B.\ Rester} \author[USM]{\quad J.\ V.\ Lambers} 
 \address[USM]{School of Mathematics and Natural Sciences, The University of Southern
Mississippi, Hattiesburg, MS 39406, USA}

\author[NRL]{\quad P.\ M.\ Jordan\corref{cor1}}
\cortext[cor1]{Corresponding author.} 

\address[NRL]{Acoustics Division, U.S.\ Naval Research Laboratory, Stennis Space Center, MS 39529, USA}


\date{\today}

\begin{abstract}
We investigate the propagation of acoustic singular surfaces, specifically, linear shock waves and nonlinear acceleration waves, in a class of inhomogeneous gases whose ambient mass density varies exponentially.  Employing the mathematical tools of singular surface theory, we first determine the evolution of both the jump amplitudes and the locations/velocities of their associated wave-fronts, along with a variety of related analytical results.  We then turn to what have become known as Krylov subspace spectral (KSS) methods to numerically simulate the evolution of the full waveforms under consideration. These are not only performed quite efficiently, since KSS allows the use of `large' CFL numbers, but also quite accurately, in the sense of capturing theoretically-predicted features of the solution profiles more faithfully than other time-stepping methods, since KSS customizes the computation of the components of the solution corresponding to the different frequencies involved.  The presentation concludes with a listing of possible, acoustics-related, follow-on studies.

 \end{abstract}

\maketitle



\section{Introduction}\label{sect:Intro}


For the general case of an inhomogeneous fluid undergoing isentropic\footnote{That is, $D\eta/Dt = 0$, where $\eta$ is the specific entropy; see, e.g., Thompson~\citep[p.~60]{T72}.} flow,  the  (compressible)  Euler system assumes the form~\cite{Bergmann46}:
\begin{subequations}\label{sys:Comp}
\begin{align}
D \rho/Dt &= -\rho (\boldsymbol{\nabla} \boldsymbol{\cdot} \textbf{v}),\label{eq:cont_1}\\
\rho D \textbf{v}/Dt  &= -\nabla p + \rho {\bf B}, \label{eq:mom_1}\\
D p/Dt &= c^2 D\rho/Dt.  \label{eq:EoS_1}
\end{align}
\end{subequations}
Here, in Cartesian coordinates, $\textbf{v}=(\mathfrak{u},\upsilon,w)$ is the velocity vector; $\rho (>0)$ is the mass density; $p (>0)$ is the thermodynamic pressure; ${\bf B}={\bf B}(x,y,z)$ is the external (per unit mass) body force  vector;    $D/Dt$ is the material derivative; and the (thermodynamic) variable $c(>0)$ denotes the sound speed.  If the fluid in question is a gas, specifically, a \emph{perfect gas}, then $c^2 = \gamma p/\rho$ and $p$, $\rho$, and $\theta$ satisfy  the following special case of the \emph{ideal gas law}~\cite[\S\,2.5]{T72}:
\be\label{eq:EoS_ideal}
p=(c_p-c_v)\rho \theta \qquad (c_{p}, c_{v} := \textrm{const.}).
\en
Here, $\theta(>0)$ is the absolute temperature; $c_p > c_v >0$ are the specific heats at constant pressure and volume, respectively; and $\gamma = c_p/c_v$,  where  $\gamma \in (1, 5/3]$ for perfect gases. 

To simplify the analyses to be performed, we  hereafter limit our attention to propagation in (perfect) gases and, moreover, always take the ambient state of the gas in question   to be \emph{quiescent}~\cite[p.~14]{P89}; i.e., while $p_{\rm a}$,  $\rho_{\rm a}$,   $\theta_{\rm a}$, and $\eta_{\rm a}$ may vary with, at most, position,  $\textbf{v}_{\rm a}  = (0,0,0)$, where  a subscript `a'   denotes the ambient state value of the quantity to which it is attached.  Our focus will be on the propagation of (acoustic) \emph{singular surfaces}~\cite{TT60}, specifically, shock and acceleration waves, in gases whose ambient density (i.e., $\rho_{\rm a}$) varies exponentially; see, e.g., Refs.~\cite{JAMM67,KJC18}, \cite[\S\,309]{Lamb16}, and~\cite{Walsh69}, as well as those cited therein,  for atmospheric propagation studies in which  this form of inhomogeneity arises.  What is more, we shall  apply our analytical and numerical tools to only the  1D versions of the acoustic models considered; specifically, those that describe (1D)  propagation along the  $z$-axis,  in Section~\ref{sect:Shock}, and the  $x$-axis, in Section~\ref{sect:ihLWE}. 

The simulation of these singular surface phenomena through numerical methods involves the use
of time-stepping schemes for partial differential equations (PDEs) 
that are required to use very small time step sizes, due to
the coupling of low- and high-frequency components of solutions---a problem described
as {\em stiffness}.  For modeling  acoustic singular surfaces in particular, high spatial 
resolution is required to represent these waves.  This, of course, stems from the fact that all singular surface 
waveforms   are of zero `thickness', thus increasing the stiffness and shrinking the time
steps even further, so that the CFL condition---a relationship between the time step size,
the spatial mesh size, and the wave speed---is satisfied.  This confluence of factors results in substantial
computational expense, and this problem is exacerbated by increasing spatial resolution.

This dilemma can be addressed by customizing the computation of components of 
the solution corresponding to different frequencies, which has led to the development
of Krylov subspace spectral (KSS) methods~\cite{block,waveblock}.  KSS methods are
designed to bring to variable-coefficient PDEs, defined on rectangular domains, the 
advantages of Fourier spectral methods for their constant-coefficient counterparts.  For
problems in the latter category, Fourier spectral methods effectively decouple components
corresponding to different frequencies, thus allowing for efficient, scalable computation 
of solutions, at least for problems with smooth solutions.  For variable-coefficient problems,
KSS methods provide similar scalability~\cite{kssepi} 
through individualized treatment of different frequencies, 
even though an explicit decoupling is not possible.  
In Ref.~\cite{nmpde}, KSS methods were applied to
1-D linear wave equations, both undamped and damped, featuring periodic media and solutions that exhibited poroacoustic shocks.  In this paper, the approach
from Ref.~\cite{nmpde} is improved and generalized to handle 1-D linear or nonlinear wave equations for modeling acoustic singular surfaces.  

In Section~\ref{sect:Shock} of this paper we consider a hybrid initial-boundary value problem (hIBVP) that models an acoustic shock in an isothermal atmosphere under a linearized special case of Sys.~\eqref{sys:Comp}.  Although the equation of motion (EoM) that results, i.e., the 1D advected damped wave equation, has constant coefficients, the jump discontinuity
that is the shock-front makes solving it via Fourier spectral methods impractical.  However, it will be demonstrated that
a KSS method, or more precisely, a Krylov subspace {\em pseudospectral} method, can solve
the PDE without the need to honor the usual CFL condition, by working in both physical and Fourier
space as needed; this is described in Section~\ref{sect:shock_scheme}.  This result is made possible by transforming the original hIBVP so that
it is feasible to represent its solution by a Fourier sine series.  In addition, the 
left boundary condition is homogenized in a manner that, in a sense, 
minimizes the effect of a newly introduced source term on the solution.

Then, in Section~\ref{sect:ihLWE}, we study  the  evolution of  acoustic acceleration waves  as they tend to finite-time blow-up, i.e., finite-time shock formation, under an inhomogeneous version of the  lossless Lighthill--Westervelt equation, a weakly-nonlinear model derivable from Sys.~\eqref{sys:Comp}.  Although KSS methods were designed for linear PDEs, it will be shown that
a KSS method can successfully be applied to this nonlinear PDE,  again without the need 
to honor the usual CFL condition.  While the combination of inhomogeneity and nonlinearity naturally poses additional challenges, as compared with the above-mentioned shock wave problem,
it will be seen that KSS methods can be effectively applied to both problems within 
a common framework that can be extended to other wave propagation problems.

In Section~\ref{sect:Numerical}, KSS based schemes are constructed for first the (linear) shock case and then the (nonlinear) acceleration wave case; see Sections~\ref{sect:shock_scheme}   
and~\ref{sect:LWE_scheme}, respectively.  This section then concludes with the presentation, in  the context of the respective propagation problems described above, of numerical simulations depicting snapshots in the evolution of the solution profiles exhibiting the former and latter singular surfaces. 
It will then be demonstrated that other time-stepping methods have difficulty resolving these singular
surfaces.

And lastly,  in Section~\ref{sect:Close}, we  record final remarks/observations and note possible follow-on investigations.




\section{Linear shock waves in an isothermal atmosphere with Rayleigh dissipation}\label{sect:Shock}

\subsection{Governing system with Rayleigh dissipation}

On setting ${\bf B}={\bf B}_{1}-\hat{\mu}\textbf{v}$ and then  linearizing about the ambient state, Sys.~\eqref{sys:Comp} becomes
\begin{subequations}\label{sys:RL-3D_1st}
\begin{align}
\partial \varrho/\partial t  + \boldsymbol{\nabla} \boldsymbol{\cdot} (\rho_{\rm a}\textbf{v})  &= 0,\label{eq:cont_3D_u}\\
 \rho_{\rm a}\partial \textbf{v}/\partial t +\boldsymbol{\nabla}\mathsf{P}  & = - [\boldsymbol{\nabla} p_{\rm a} - (\varrho +\rho_{\rm a}){\bf B}_{1}+\hat{\mu}\rho_{\rm a}\textbf{v}], \label{eq:mom_3D_u}\\
\partial \mathsf{P}/\partial t +  \textbf{v} \boldsymbol{\cdot} (\boldsymbol{\nabla} p_{\rm a}) &= c_{\rm a}^{2} [\partial \varrho/\partial t +  \textbf{v} \boldsymbol{\cdot} (\boldsymbol{\nabla} \rho_{\rm a})]. \label{eq:Egy_3D_u}
\end{align}
\end{subequations}
Here, $\mathsf{P} = p - p_{\rm a}$ and $\varrho =\rho - \rho_{\rm a}$, while the gas law and sound speed expressions  become
\be\label{eq:EoS_perfect}
  p_{\rm a}= (c_p-c_v)\rho_{\rm a}\theta_{\rm a} \qquad (c_{p}, c_{v} := \text{const.}),
\en
\be\label{eq:c_ambient}
c_{\rm a}=\sqrt{\gamma p_{\rm a}/ \rho_{\rm a}}= \sqrt{c_{p}(\gamma-1)\theta_{\rm a}},
\en
respectively.  In writing Sys.~\eqref{sys:RL-3D_1st} we have  introduced \emph{Rayleigh's  resistance law}\footnote{See, e.g., Lamb~\cite[p.~389]{Lamb16}, who credits this `artifice' to Lord Rayleigh and notes, among other things, its usefulness in approximating `small dissipative forces'.  Also, Rayleigh's resistance law  should not to be confused with any of the numerically-motivated  \emph{artificial  viscosity} methods~\cite{Roache72} put forth since 1950.}, i.e., the term $-\hat{\mu}\textbf{v}$, so as to prevent the shock amplitude blow-up that would otherwise afflict the solution of the hIBVP considered in the next subsection; see, e.g., Ref.~\cite[\S\,5.4]{PMJ22}.  Here, the  (resistance) coefficient $\hat{\mu}(>0)$ carries units of 1/s.

Sys.~\eqref{sys:RL-3D_1st} can be reduced to a two-equation system as follows: First make use of the equilibrium condition
\be\label{eq:grad_pa}
\boldsymbol{\nabla} p_{\rm a}=\rho_{\rm a}{\bf B}_{1},
\en
by which we  eliminate $\boldsymbol{\nabla} p_{\rm a}$, and then use  Eq.~\eqref{eq:cont_3D_u} to eliminate $\partial \varrho/\partial t$ from Eqs.~\eqref{eq:mom_3D_u} and~\eqref{eq:Egy_3D_u}, after applying $\partial/\partial t$ to the former; these actions yield, after simplifying, 
\begin{subequations}\label{sys:RL-3D_2eqs}
\begin{align}
 \rho_{\rm a}\partial^2 \textbf{v}/\partial t^2 +\boldsymbol{\nabla}(\partial \mathsf{P}/\partial t)  & = - [\boldsymbol{\nabla} \boldsymbol{\cdot} (\rho_{\rm a} \textbf{v})]{\bf B}_{1} - \hat{\mu}\rho_{\rm a}\partial \textbf{v}/\partial t, \label{eq:mom_3D_RL-2}\\
\partial \mathsf{P}/\partial t &= -  \textbf{v}\boldsymbol{\cdot} (\rho_{\rm a}{\bf B}_{1})  - \rho_{\rm a}c_{\rm a}^{2} (\boldsymbol{\nabla} \boldsymbol{\cdot} \textbf{v}). \label{eq:Egy_3D_RL-2}
\end{align}
\end{subequations}

The final step is, of course, eliminating $\partial \mathsf{P}/\partial t$ between the  PDEs that comprise Sys~\eqref{sys:RL-3D_2eqs}; as this is easily accomplished, we omit the details and present the  (vector) EoM that results, viz.:
\be\label{eq:vecEoM_shock}
\rho_{\rm a}\partial^2 \textbf{v}/\partial t^2 -\gamma \boldsymbol{\nabla}(p_{\rm a} \boldsymbol{\nabla} \boldsymbol{\cdot} \textbf{v}) +  \hat{\mu}\rho_{\rm a}\partial  \textbf{v}/\partial t = \boldsymbol{\nabla}[ \textbf{v} \boldsymbol{\cdot} (\rho_{\rm a}{\bf B}_{1})] - [\boldsymbol{\nabla} \boldsymbol{\cdot} (\rho_{\rm a} \textbf{v})]{\bf B}_{1};
\en
compare  with the   EoM derived in  Ref.~\cite[\S\, 311]{Lamb16}, which is stated in component form, and wherein each component of the body force vector (i.e., $(X, Y, Z)$) is assumed constant.

\subsection{Vertical propagation in a 1D, isothermal atmosphere}\label{sect:Isothermal}

This case is defined by the following assumptions: 
\begin{itemize}
\item[(I)] $\theta_{\rm a}(z) =\theta_{0}$, which via  Eq.~\eqref{eq:c_ambient} yields $c_{\rm a}(z)=c_{0}$, where 
\be\label{eq:c0}
c_{0}= \sqrt{c_{p}(\gamma-1)\theta_{0}}.
\en
\item[(II)]  ${\bf B}_{1}=(0,0,-g)$, where $g$ represents the magnitude of the (constant) acceleration due to gravity near the surface; in the case of Earth,  $g \approx 9.81\textrm{m/s}^2$.
\item[(III)]   1D  propagation along the $z$-axis; i.e., the velocity vector simplifies to $\textbf{v} =(0, 0, w(z,t))$,  while  $p=p(z, t)$, $\rho=\rho(z, t)$,  and $\theta = \theta(z, t)$, where the $+z$-axis is taken to be directed vertically upward. Note also that  $p_{\rm a}$, $\rho_{\rm a}$, $\vartheta_{\rm a}$, and $c_{\rm a}$ are now (at most) functions of $z$ only.
\item[(IV)] The `flat Earth' approximation (see, e.g., Ref.~\cite{LH98}) and the negligibility of the Earth's rotation.
\end{itemize}
Here,  $\theta_{0}=(p_{0}/\rho_{0})/(c_p-c_v)$, where $p_{0}, \rho_{0} = \lim_{z \to 0}p_{\rm a}(z), \rho_{\rm a}(z)$, respectively, are both positive.

Under the first three of these assumptions, Eq.~\eqref{eq:grad_pa} reduces to 
\be
\frac{\rd \rho_{\rm a}}{\rd z}=-\left[\frac{g}{c_{v}(\gamma -1)\theta_{0}} \right]\rho_{\rm a};
\en
integrating this simple ODE subject to $\rho_{\rm a}(0)=\rho_{0}$ yields 
\be\label{eq:rho_a_Iso}
\rho_{\rm a}(z) =\rho_{0}\exp(-z/H),
\en 
and this in turn yields,  by way Eq.~\eqref{eq:EoS_perfect} and our assumption $\theta_{\rm a}(z) =\theta_{0}$,
\be
 p_{\rm a}(z) =p_{0}\exp(-z/H).
 \en
Here,  following convention, we have set $H=c_{0}^{2}/(\gamma g)$, where $H$ is known  as the `scale height of the atmosphere'~\cite[p.~69]{T72}.  



\subsection{Problem formulation and solution}

In this subsection we consider  the following \emph{hybrid}\footnote{See Ref.~\cite[Footnote~4]{KJC18}.} initial-boundary value problem~(hIBVP), whose EoM is the result of applying the above-listed assumptions/consequences to  Eq.~\eqref{eq:vecEoM_shock}:
\begin{subequations}\label{IBVP:Atmos_Isothermal}
\begin{align}
\label{eq:Atmos_Ray_EoM}
&  w_{tt}  - c_{0}^{2}w_{zz}  + \gamma g w_{z} +\hat{\mu} w_{t} =   0, \quad (z,t)\in (0, \infty)\times (0, \infty),\\
& w(0,t) = W_{0}\mathcal{H} (t)\cos(\omega t),\quad \lim_{z \to \infty}\overline{w}(z, s) =0, \quad  t, s >0,\\
& w(z,0) = 0, \quad  w_{t}(z,0) = 0, \quad  z >0.
\end{align}
\end{subequations}
In hIBVP~\eqref{IBVP:Atmos_Isothermal}, which is commonly referred to in the acoustics literature as a `signaling problem'; $\mathcal{H}(\cdot)$ is the Heaviside unit step function;  the (known) constants $W_{0}(> 0)$ and $\omega(>0)$ are the amplitude and angular frequency, respectively, of the inserted signal; and, in the present section, $\overline{w}(z, s)$ denotes the image of $w(z,t)$ in the Laplace  domain, i.e.,
\be\label{eq:LapTrans}
\overline{w}(z, s) =\int_{0}^{\infty}\exp(-s t) w(z,t)\,\rd t \qquad (s>0),
\en
where in this section $s$ denotes the Laplace transform parameter; see, e.g., Ref.~\cite{CJ63}.


To determine the evolution of the shock wave modeled by this problem, we begin by applying  the Laplace transform, i.e., Eq.~\eqref{eq:LapTrans},  to  Eq.~\eqref{eq:Atmos_Ray_EoM}, which we observe is the (1D) advected damped wave equation~\cite{Mickens22}, and the boundary condition (BC) at $z=0$.  After employing the initial conditions  and simplifying,  we are led to consider the subsidiary equation
\be
c_{0}^{2} \overline{w}^{\prime\prime}-\gamma g  \overline{w}^{\prime}-(s^2+\hat{\mu} s)\overline{w}=0,
\en
where in this section a prime denotes $\rd/\rd z$.  Solving this ODE subject to the transformed BC, which reads $\overline{w}(0,s) = W_{0}s/(s^2+\omega^2)$, and the asymptotic condition in the Laplace domain, i.e., $\lim_{z \to \infty}\overline{w}(z, s)=0$,  yields, after some manipulation, 
\be\label{eq:gas_Ray_LapDomian}
\overline{w}(z,s)=\left(\frac{W_{0}s}{s^2+\omega^2}\right)\exp\left[z/(2H)\right]\exp\!\left[-\left(\frac{z}{c_{0}}\right)\sqrt{\left(s+\frac{1}{2}\hat{\mu} \right)^2-\chi^2}\, \right]\!,
\en
which of course is the exact solution of hIBVP~\eqref{IBVP:Atmos_Isothermal} in the Laplace  domain; here, on the assumption that $\hat{\mu} > c_{0}/H$ (see Section~\ref{sect:shock_small-t} below), we have set
\be
\chi := \frac{1}{2}\sqrt{\hat{\mu}^2-\frac{c_{0}^{2}}{H^{2}}}.
\en

Due to the relative simplicity of Eq.~\eqref{eq:gas_Ray_LapDomian}, its exact inverse can be obtained  in several different ways; e.g.,  using a suitably extensive table of Laplace inverses  and the properties of the Laplace transform.  Here,  we employ the results given in  Ref.~\cite[\S\,90]{CJ63} in conjunction with the convolution theorem for Laplace transforms to show that
\begin{multline}\label{eq:Ray_t-domain-sol}
w(z,t) = W_{0}\mathcal{H}(t - z/c_{0})\Bigg\{ \exp\left[-\left(\frac{\hat{\mu} - c_{0}/H}{2}\right)\frac{z}{c_{0}}\right] \cos[\omega(t-z/c_{0})]\\
+\frac{\chi z}{c_{0}}\int_{z/c_{0}}^{t}\frac{\exp(-\hat{\mu} \varsigma/2)\cos[\omega (t-\varsigma)]I_{1}\left[\chi\sqrt{\varsigma^2-(z/c_{0})^2}\,\right]}{\sqrt{\varsigma^2-(z/c_{0})^2}}\,\rd \varsigma \Bigg\},
\end{multline}
 where in this communication $I_{\nu}(\cdot)$ denotes the modified Bessel function of the first kind of order $\nu$. 
 
 \subsection{Shock analysis and small-time results}\label{sect:shock_small-t}
 
An inspection of Eq.~\eqref{eq:Ray_t-domain-sol} reveals that $w$ exhibits a jump discontinuity of amplitude
 \be\label{eq:shock-Amp_Isothermal}
 \lshad w\rshad(t) =W_{0} \exp\left[-\left(\frac{\hat{\mu} - \hat{\mu}_{\rm c}}{2}\right)t\right],
 \en
 where $\lshad \cdot\rshad(t)$ is defined in Appendix~\ref{App:SST} and  $\hat{\mu}_{\rm c} := c_{0}/H=\gamma g/c_{0}$ is a critical value of $\hat{\mu}$, across the 
 \emph{shock wave}\footnote{As defined in Appendix~\ref{App:SST}.}  $z=\Sigma(t)$, where in this section
 \be\label{eq:shock-Loc_Isothermal}
 \Sigma(t) =c_{0}t;
 \en
 here, we observe that $\Sigma(t)$ propagates to the right along the $z$-axis with speed
 \be
 \frac{\rd \Sigma(t)}{\rd t}=c_{0}.
 \en
From Eq.~\eqref{eq:shock-Amp_Isothermal} it is clear that we must take $\hat{\mu} > \hat{\mu}_{\rm c}$ if we are to have $\lshad w\rshad(t) \to 0$  as  $t\to \infty$. Eq.~\eqref{eq:shock-Amp_Isothermal} also makes clear that $\lshad w\rshad(t)$ is independent of $\omega$; see Eq.~\eqref{eq:Ray_small-t} below.

 The  transient nature of $\lshad w\rshad(t)$, which follows from our assumption $\hat{\mu} > \hat{\mu}_{\rm c}$, and the complicated nature of Eq.~\eqref{eq:Ray_t-domain-sol}  motivate us to seek a simpler, approximate expression for the latter---one that is valid for small-$t$. To this end, we turn to Ref.~\cite[\S\S\,123, 124]{CJ63} and expand Eq.~\eqref{eq:gas_Ray_LapDomian} for large-$s$;  this  yields, after simplifying, 
\begin{multline}\label{eq:Ray_large-s}
\overline{w}(z,s)\sim W_{0} s^{-1}\exp\left(-sz/c_{0}\right)\exp\left[-\left(\hat{\mu} - \frac{c_{0}}{H}\right)\frac{z}{2c_{0}} \right] \Bigg{\{} 1+ \left( \frac{\hat{\mu}^2-c_{0}^{2}/H^2}{8c_{0}}\right)\!\frac{z}{s} \\
+ \left[ \frac{1}{2}\left(\frac{\hat{\mu}^2-c_{0}^{2}/H^2}{8}\right)^{2}\left(\frac{z}{c_{0}}\right)^{2}- \frac{\hat{\mu}}{2} \left(\frac{\hat{\mu}^2-c_{0}^{2}/H^2}{8}\right)\left(\frac{z}{c_{0}}\right)-\omega^{2}\right]\!\frac{1}{s^2}
+\mathcal{O}(s^{-3})\Bigg{\}} \qquad (s \to \infty).
\end{multline}
Inverting term-by-term using a standard table of inverses ~\cite{CJ63}, we find that
\begin{multline}\label{eq:Ray_small-t}
w(z,t) \approx  W_{0}\mathcal{H}(t - z/c_{0})\exp\left[-\left(\frac{\hat{\mu} -\hat{\mu}_{\rm c}}{2}\right)\frac{z}{c_{0}}\right]\Bigg{\{} 1+\left( \frac{\hat{\mu}^2- \hat{\mu}_{\rm c}^2}{8}\right)\left(\frac{z}{c_{0}}\right) (t - z/c_{0}) \\
+\frac{1}{2}\left[ \frac{1}{2}\left(\frac{\hat{\mu}^2- \hat{\mu}_{\rm c}^{2}}{8}\right)^{2}\!\left(\frac{z}{c_{0}}\right)^{2}- \frac{\hat{\mu}}{2} \left(\frac{\hat{\mu}^2- \hat{\mu}_{\rm c}^{2}}{8}\right)\left(\frac{z}{c_{0}}\right)-\omega^{2}\right] (t - z/c_{0})^{2} \Bigg{\}},
\end{multline}
where terms of $\mathcal{O}[(t - z/c_{0})^{3}]$ have been neglected.  It should be noted that while Eq.~\eqref{eq:Ray_small-t} can be expected to prove accurate at, and in the  nearby-vicinity of, $z=\Sigma (t)$ for all $t>0$, its usefulness as a global approximation to $w$ is limited to sufficiently small values of $t$.

We conclude this subsection with the following observations: (i) Note that  $\hat{\mu}_{\rm c}=1/t_{H}$, where $\Sigma(t_{H})=H$; i.e., the critical value of $\hat{\mu}$ is the reciprocal of the time taken by $\Sigma$ to propagate from $z=0$ to $z=H$. (ii) For $t$ sufficiently small, the impact of $\omega$ is of only second order. And (iii), Eq.~\eqref{eq:shock-Amp_Isothermal} could have also been obtained directly from Eq.~\eqref{eq:Ray_large-s} via the application of the theorem in Ref.~\cite[\S\,4]{BH68}.  


\subsection{Large-time asymptotic results}\label{sect:Large-t}

Using  the properties of the Laplace transform detailed in Ref.~\cite[\S\,126]{CJ63}, it is not difficult to establish that, to lowest order, the large-$t$ behavior  of Eq.~\eqref{eq:Ray_t-domain-sol} is described by
\be\label{eq:Large-t}
w(z,t) \sim  W_{0}\exp\left[-\left(\frac{\sigma-\hat{\mu}_{\rm c}}{2}\right)\frac{z}{c_{0}} \right] \cos\left[\omega t- \tfrac{1}{2}(\varkappa/c_{0})z\right] \qquad    (t \to \infty),
\en
where
\be
\sigma = \sqrt{\frac{\hat{\mu}_{\rm c}^{2} - 4\omega^2+\sqrt{(\hat{\mu}_{\rm c}^{2} - 4\omega^2)^2 +16(\hat{\mu}\omega)^{2}}}{2}} \quad \text{and} \quad \varkappa =\sqrt{\frac{-\hat{\mu}_{\rm c}^{2} + 4\omega^2+\sqrt{(\hat{\mu}_{\rm c}^{2} - 4\omega^2)^2 +16(\hat{\mu}\omega)^{2}}}{2}}.
\en
It can be shown that $\sigma \in (\hat{\mu}_{\rm c}, \hat{\mu})$ for all $\omega >0$, where we recall our assumption $\hat{\mu}_{\rm c} < \hat{\mu}$; thus, provided $\hat{\mu}_{\rm c} < \hat{\mu}$ holds, $w(z,t) \to 0$, as $z \to \infty$, under Eq.~\eqref{eq:Ray_t-domain-sol}. 

Also noteworthy is the fact that Eq.~\eqref{eq:Large-t} admits the high-frequency simplification:
\be\label{eq:Large-t_omega}
w(z,t) \sim  W_{0}\exp\left[-\left(\frac{\hat{\mu} - \hat{\mu}_{\rm c}}{2}+\mathcal{O}(\omega^{-2})\right)\frac{z}{c_{0}} \right] \cos\left[\omega t- \omega\left(1+\frac{\hat{\mu}^{2} - \hat{\mu}_{\rm c}^{2}}{8\omega^{2}}+\mathcal{O}(\omega^{-4})\right) \frac{z}{c_{0}} \right] \qquad    (t \to \infty),
\en
a waveform with phase velocity
\be\label{eq:Vph}
V_{\rm ph} \sim c_{0}\left[1-\frac{\hat{\mu}^{2} - \hat{\mu}_{\rm c}^{2}}{8\omega^{2}}+\mathcal{O}(\omega^{-4})\right]  \qquad    (\omega \to \infty). 
\en
Eq.~\eqref{eq:Vph} indicates that, although attenuated, our waveform exhibits  \emph{dispersion-free} propagation in the limit $\omega \to \infty$.

We conclude this subsection by asking the reader to compare Eq.~\eqref{eq:Large-t_omega} with the large-$t$ expression derived in Ref.~\cite[(9) \S\,4$\cdot$17]{McL62} for a variant of hIBVP~\eqref{IBVP:Atmos_Isothermal}; therein, among other differences,  `$\beta$' plays the role of  $\tfrac{1}{2}\hat{\mu}$, where `$\beta$' should be regarded as though it were  independent of `$\omega$', and $\hat{\mu}_{\rm c} :=0$.
  


\section{Nonlinear effects: Acceleration waves under a variable-coefficient Lighthill--Westervelt equation}\label{sect:ihLWE}
 
Assuming $p_{\rm a}({\bf x})\equiv p_{0}$, the lossless version of the inhomogeneous Lighthill--Westervelt equation (ihLWE) can be expressed as~\cite[Eq.~(2.55)]{Reiso91}
\be\label{eq:LWE_wp_3D}
 \wp_{tt} - c_{\rm a}^{2}({\bf x})\nabla^{2}\wp +  \rho_{\rm a}^{-1}({\bf x})c_{\rm a}^{2}({\bf x}) \left\{[\boldsymbol{\nabla} \rho_{\rm a}({\bf x})] \boldsymbol{\cdot \nabla}\wp\right\} 
 =\beta \rho_{\rm a}^{-1}({\bf x}) c_{\rm a}^{-2}({\bf x})\partial_{tt}(\wp^{2});
\en
see also Ref.~\cite[Eq.~(2)]{JASA01}, taking note of the fact that `$\kappa_{0}$' therein is equal to $1/(\gamma p_{0})$ in the present section.  In Eq.~\eqref{eq:LWE_wp_3D},  $\wp({\bf x}, t) = p({\bf x}, t) - p_{0}$  denotes the acoustic pressure; $\beta$, the coefficient of nonlinearity, is given by $\beta =(\gamma+1)/2$ in the case of perfect gases;  and in this communication we let 
\be
\partial_{\varsigma}, \partial_{\varsigma \varsigma} := \partial^{1,2}/\partial \varsigma^{1,2},
\en
respectively.  This weakly-nonlinear EoM is readily derived from 
Sys.~\eqref{sys:Comp} after setting ${\bf B}=(0,0,0)$, which  follows from our assumptions $p_{\rm a}({\bf x})\equiv p_{0}$ and $\textbf{v}_{\rm a}=(0,0,0)$, and invoking the \emph{finite-amplitude approximation scheme} (FAAS); see, e.g., Ref.~\cite[\S\, 2.1]{Reiso91}\footnote{The system stated in Ref.~\cite[p.~3]{Reiso91} reduces to a form equivalent to Sys.~\eqref{sys:Comp} on setting $\kappa, \mu, K, q =0$, and replacing ${\bf F}$ with $\rho {\bf B}$, in the former.}.  In carrying out the analysis in this section, we will make use of the fact that Eq.~\eqref{eq:LWE_wp_3D} can also be expressed as
\be\label{eq:LWE_Psi_3D}
c_{\rm a}^{2}({\bf x})\nabla^{2}\Phi - \Phi_{tt}+  c_{\rm a}^{2}({\bf x}) \boldsymbol{\nabla \cdot} 
\left \{[\rho_{\rm a}^{-1}({\bf x}) \boldsymbol{\nabla} \rho_{\rm a}({\bf x})] \Phi \right \} =\beta c_{\rm a}^{-2}({\bf x})\partial_{t}[(\Phi_{t})^{2}];
\en
here, $\textbf{v}=\nabla \Phi$, where $\Phi = \Phi({\bf x},t)$ is the scalar velocity potential, and we note that $\wp$, $\Phi$ are related via
\be\label{eq:Phi_wp}
\wp = -\rho_{\rm a}({\bf x})\Phi_{t} +\text{h.o.t.}
\en

In the case of 1D propagation along the $x$-axis, $\textbf{v}=(\mathfrak{u}(x,t),0,0)$, where now $\mathfrak{u}(x,t)=\Phi_{x}(x,t)$, and the acoustic pressure simplifies to $\wp(x, t) = p(x, t) - p_{0}$, while $\rho_{\rm a}$ and $c_{\rm a}$ become functions of (at most) $x$; consequently, Eqs.~\eqref{eq:LWE_wp_3D} and~\eqref{eq:LWE_Psi_3D} are reduced to
\be\label{eq:LWE_wp_1D}
\wp_{tt} - c_{\rm a}^{2}(x)\wp_{xx}  +  c_{\rm a}^{2}(x) [\rho_{\rm a}^{\prime}(x)/\rho_{\rm a}(x)]\wp_{x} = \beta \rho_{\rm a}^{-1}(x) c_{\rm a}^{-2}(x)\partial_{tt}(\wp^{2}),
\en
\be\label{eq:LWE_Psi_1D}
c_{\rm a}^{2}(x)\Phi_{xx} - \Phi_{tt}+  c_{\rm a}^{2}(x) \partial_{x} \left\{[\rho_{\rm a}^{\prime}(x)/\rho_{\rm a}(x)] \Phi \right \} = \beta c_{\rm a}^{-2}(x)\partial_{t}[(\Phi_{t})^{2}],
\en
respectively, where in this section a prime denotes $\rd/\rd x$.

If we now assume the ambient density of the gas exhibits the following \emph{exponential} profile:
\be
\rho_{\rm a}(x)=\rho_{0}\exp\left(\alpha x/\ell \right), 
\en
then, by Eq.~\eqref{eq:c_ambient} and our earlier assumption $p_{\rm a}(x)=p_{0}$, if follows that $c_{\rm a}(x)=\hat{c}_{\rm a}(x)$, where
\be
\hat{c}_{\rm a}(x)= c_{0}\exp\left[-\alpha x/(2\ell)\right].
\en
Here, $c_{0}=\sqrt{\gamma p_{0}/\rho_{0}}$, as in the previous section;  $\alpha(\neq  0)$ is a free (dimensionless) parameter;  and $\ell(>0)$ is the length of the propagation domain.  On making these substitutions  and simplifying, Eqs.~\eqref{eq:LWE_wp_1D} and~\eqref{eq:LWE_Psi_1D} reduce to
\be\label{eq:LWE_p}
\wp_{tt} - \hat{c}_{\rm a}^{2}(x)\wp_{xx}  + (\alpha/\ell) \hat{c}_{\rm a}^{2}(x)\wp_{x}= p_{0}^{-1}(\beta/\gamma)
\partial_{tt}(\wp^{2}),
\en
\be\label{eq:LWE_Psi}
\hat{c}_{\rm a}^{2}(x)\Phi_{xx} - \Phi_{tt}+  (\alpha/\ell) \hat{c}_{\rm a}^{2}(x) \Phi_{x} = \beta \hat{c}_{\rm a}^{-2}(x)\partial_{t}[(\Phi_{t})^{2}],
\en
respectively.

We seek to analytically investigate, and later numerically solve, Eq.~\eqref{eq:LWE_p} subject to the following BCs and ICs:
\begin{multline}
\wp(0,t)=p_{\rm pk}\mathcal{H} (t)\sin(\omega t), \quad \wp(\ell,t)=0 \qquad (0<t< t_{\rm f}); \qquad
\wp(x,0)=0, \quad \wp_{t}(x,0)=0 \qquad (0<x< \ell),
\end{multline}
where, adopting the notation of Ref.~\cite[\S\,1-8]{P89},   the parameters $p_{\rm pk}(>0)$ and $\omega(>0)$ denote the peak pressure (or amplitude) and angular frequency, respectively, of the inserted signal.  Moreover, we let $t_{\rm f}(>0)$ represent the  final (i.e., largest) instant of time at which we seek to compute $\wp$, and we recall that $\mathcal{H} (\cdot)$ denotes the Heaviside unit step function.  
 
To this end, we   introduce the following dimensionless variables:
\be\label{eq:ND_variables}
P=\wp/p_{\rm pk}, \qquad X=x/\ell, \qquad T = (c_{0}/\ell)t, \qquad \mathsf{U} = \mathfrak{u}/\mathfrak{u}_{\rm c},
\en
where $\mathfrak{u}_{\rm c}(>0)$ denotes a characteristic value of the velocity field,  and  consider the following (dimensionless) signaling IBVP involving the ihLWE:
\begin{subequations}\label{IBVP:LWE_non-D}
\begin{align}
\label{eq:EoM_P}
&\left(1- 2\epsilon\hat{\beta} P \right) \! P_{TT} - V_{\rm a}^{2}(X)P_{XX}  +\alpha V_{\rm a}^{2}(X)P_{X}  = 2\epsilon\hat{\beta} (P_{T})^{2},  \qquad  (X,T)\in (0, 1)\times (0, T_{\rm f}),\\
&P(0,T) = \mathcal{H} (T)\sin(\pi T),\quad P(1,T) =0, \qquad  T\in (0, T_{\rm f}), \\
&P(X,0) = 0, \quad  P_{T}(X,0) = 0, \qquad  X\in (0,1).
\end{align}
\end{subequations}
In IBVP~\eqref{IBVP:LWE_non-D}, $\epsilon = p_{\rm pk}/p_{0}$, where $\epsilon \ll 1$ is assumed in accordance with the FAAS; we have set 
\be\label{eq:beta-hat}
\hat{\beta} :=\beta/\gamma = (1+1/\gamma)/2,
\en
where $\hat{\beta}\in [4/5, 1)$; the dimensionless form of  $\hat{c}_{\rm a}^{2}(x)$ is
\be
V_{\rm a}^{2}(X) =\exp(-\alpha X);
\en
we have taken the dimensionless signal frequency equal to $\pi$;   $T_{\rm f}$, the dimensionless version of $t_{\rm f}$,   is given by
\be\label{eq:Tf}
T_{\rm f}=\begin{cases}
\begin{cases}
T_{1}, & \alpha \in (-\infty,0)\cup (0,\alpha^{\bullet}],\\
T_{\infty}, & \alpha >\alpha^{\bullet},
\end{cases} & \epsilon \in (0, \epsilon^{\bullet}),\\
\begin{cases}
T_{1}, & \alpha < 0,\\
T_{\infty}, & \alpha >0,
\end{cases} & \epsilon =\epsilon^{\bullet},\\
\begin{cases}
T_{1}, & \alpha < -|\alpha^{\bullet}|,\\
T_{\infty}, & \alpha \in [-|\alpha^{\bullet}|, 0)\cup (0, +\infty),\\
\end{cases} & \epsilon > \epsilon^{\bullet},
\end{cases}
\qquad \alpha \neq 0,
\en
where we have set 
\be\label{eq:epsilon_bullet}
\epsilon^{\bullet} := 1/(\pi\hat{\beta});
\en
 and  $T_{1}$, $T_{\infty}$, and $\alpha^{\bullet}$ are defined below in Eqs.~\eqref{eq:T1}, \eqref{eq:Tinfty}, and~\eqref{eq:alpha_bullet}, respectively.

Now making use of Eq.~\eqref{eq:Phi_wp} and the relation $\mathfrak{u}(x,t)=\Phi_{x}(x,t)$, it is not difficult to decompose Eq.~\eqref{eq:LWE_Psi} into the (kinematic-like) system
\begin{subequations}\label{sys:ihLWE}
\begin{align}
\epsilon\left(1-2\epsilon \hat{\beta}P\right)P_{T}+\gamma \textrm{Ma}\mathsf{U}_{X} &=-(\gamma \alpha\textrm{Ma})\mathsf{U},\label{eq:ihLWE1}\\
\gamma\textrm{Ma}\exp(\alpha X)\mathsf{U}_{T} + \epsilon P_{X} &= \epsilon\alpha P, \label{eq:ihLWE2}
\end{align}
\end{subequations}
which we have expressed in terms of the dimensionless variables defined in Eq.~\eqref{eq:ND_variables}. In Sys.~\eqref{sys:ihLWE},  all h.o.t.\ terms have been neglected, in accordance with the FAAS, and $\textrm{Ma} = \mathfrak{u}_{\rm c}/c_{0}$ is the Mach number, where, also based on the FAAS, $\textrm{Ma} \ll 1$ is  assumed.


In the case of IBVP~\eqref{IBVP:LWE_non-D}, it is readily established, by first taking jumps of Sys.~\eqref{sys:ihLWE}, and then applying the tools of singular surface theory (see Appendix~\ref{App:SST}) to the resulting jump equations, that the amplitudes of the jumps in the first derivatives of $P$, i.e., the acceleration wave amplitudes,  are given by 
\be\label{eq:[PT]}
\lshad P_{T} \rshad(T) = \frac{3\alpha \pi \sqrt{1+\alpha T/2}}
{(3\alpha +4\pi \epsilon\hat{\beta})-4\pi \epsilon\hat{\beta} (1+\alpha T/2)^{3/2}},
\en
\be\label{eq:[PX]}
\lshad P_{X} \rshad(T) = \frac{-3\alpha \pi (1+\alpha T/2)^{3/2}}
{ (3\alpha +4\pi \epsilon\hat{\beta}) - 4\pi \epsilon\hat{\beta} (1+\alpha T/2)^{3/2}}.
\en
Moreover, both the former and latter  occur across the the (planar) surface $X=\ram(T)$,  where 
\be\label{eq:ram}
\ram(T) =  2\alpha^{-1} \ln\left(1+ \tfrac{1}{2}\alpha  T \right).
\en
Here, $X=\ram(T)$ is the  equation of  the  \emph{acceleration wave-front}  associated with IBVP~\eqref{IBVP:LWE_non-D}; i.e., $\ram(T)$ is the dimensionless form of the expression for $\Sigma(t)$ in the present section (again, see Appendix~\ref{App:SST}). 

It should be noted that Eq.~\eqref{eq:[PT]} was determined first;  Eq.~\eqref{eq:[PX]} then followed from the former via Eq.~\eqref{eq:Max_CC}, which in the present setting assumes the form:
\be
\mathcal{V}(T)\lshad P_{X} \rshad(T) = -\lshad P_{T} \rshad(T).
\en
Here, $\mathcal{V}(T)=\rd \ram(T)/\rd T$ is easily shown to be given by
\be
\mathcal{V}(T)=\frac{1}{1+\tfrac{1}{2}\alpha T},
\en
where $\mathcal{V}(T)(>0)$ is  the velocity at which $\ram(T)$ propagates (to the right) along the $X$-axis.  Also, Eq.~\eqref{eq:ram} was obtained by integrating the following ODE subject to the initial condition $\ram(0)=0$:
\be
(\rd \ram(T)/\rd T)^{2} = V_{\rm a}^{2}(\ram(T)),
\en
where we observe that $\ram(T_{1})=1$; here, $T_{1}(>0)$, which is given by  
\be \label{eq:T1}
T_{1} := \frac{2[\exp(\alpha/2)-1]}{\alpha},
\en
 is  the (dimensionless) time at which $X=\ram(T)$  makes its \emph{initial} arrival at the right boundary (i.e., $X=1$).



Of particular interest to us, however, is the following:  An inspection of Eqs.~\eqref{eq:[PT]} and~\eqref{eq:[PX]} reveals that both amplitude expressions are capable of exhibiting finite-time shock formation; specifically, $\big{|}\lshad P_{T} \rshad(T_{\infty})\big{|} = \big{|}\lshad P_{X} \rshad(T_{\infty})\big{|}=\infty$, where the (common) blow-up time is given by
\be\label{eq:Tinfty}
T_{\infty} =-\,\frac{2}{\alpha}\left[1-\left(1+\frac{3\alpha}{4\pi \epsilon\hat{\beta}}\,\right)^{2/3}  \right].
\en 
To simplify our analysis of this phenomenon, we henceforth assume  that $\alpha = \alpha^{\bullet}$;   here, recalling Eqs.~\eqref{eq:Tf} and~\eqref{eq:epsilon_bullet}, we record the expressions
\be\label{eq:alpha_bullet}
\alpha^{\bullet} := -\,\frac{4}{3}\begin{cases}
\epsilon\hat{\beta}\pi +
\mathfrak{W}_{-1}\left[-\epsilon\hat{\beta}\pi \exp(-\epsilon\hat{\beta}\pi )\right], &  \epsilon \in (0, \epsilon^{\bullet}),\\
\\
\epsilon\hat{\beta}\pi + \mathfrak{W}_{0}\left[-\epsilon\hat{\beta}\pi \exp(-\epsilon\hat{\beta}\pi )\right], &  \epsilon > \epsilon^{\bullet},
\end{cases}
\en
which were derived by setting the right-had side of Eq.~\eqref{eq:T1} equal to the right-hand side of Eq.~\eqref{eq:Tinfty} and then solving for $\alpha$.  In Eq.~\eqref{eq:alpha_bullet},  wherein we have \emph{implicitly} imposed the restriction $\epsilon \neq \epsilon^{\bullet}$ to ensure that $\alpha^{\bullet} \neq 0$,  we have used $\mathfrak{W}_{-1,0}(\cdot)$ to denote  the `negative one' and `principal'  branches, respectively, of the Lambert $W$-function~\cite{Corless96}.

Thus, in what follows,
\be
T_{\rm f}=T_{\infty}\Big{|}_{\alpha=\alpha^{\bullet}} = T_{1}\Big{|}_{\alpha=\alpha^{\bullet}};
\en
recall Eqs.~\eqref{eq:Tf}.  A consequence of this equality is, of course, that
\be\label{eq:Bdy_shock}
\ram(T_{\infty})\Big{|}_{\alpha=\alpha^{\bullet}} =1.
\en
For $\alpha = \alpha^{\bullet}$, then, Eq.~\eqref{eq:Bdy_shock} indicates that shock formation occurs  the \emph{instant} our acceleration wave-front $X=\ram(T)$ arrives at the right boundary;  here, we observe that $\alpha^{\bullet} > 0$ (resp.~$\alpha^{\bullet} < 0$)  for $\epsilon \in (0, \epsilon^{\bullet})$ (resp.~$\epsilon > \epsilon^{\bullet}$).  

Because we have selected $\alpha = \alpha^{\bullet}$, and $\alpha^{\bullet} > \alpha_{\rm crt}$, an inequality which is easily verified, it follows that $T_{\infty}\in\mathbb{R}^{+}$; here,
\be
\alpha_{\rm crt}=-\,\frac{4\pi \epsilon\hat{\beta}}{3},
\en
where we observe that  $\alpha = \alpha_{\rm crt}$ (resp.~$\alpha < \alpha_{\rm crt}$) implies that $T_{\infty}=T_{\rm BD}$ (resp.~$T_{\infty}\in\mathbb{C}$).  What we denote by $T_{\rm BD}(>0)$ is the value of $T$ at which the \emph{breakdown} of the $\alpha <0$ case of the ihLWE occurs; it is given by
\be
T_{\rm BD} =\frac{2}{|\alpha|} \qquad (\alpha < 0),
\en
and  by `breakdown' we mean $\mathcal{V}(T), \ram(T) \to \infty$ as $T \to T_{\rm BD}$, which we stress is a possibility only when $\alpha <0$.  We hasten to point out, however, that this breakdown \emph{cannot} occur under the present formulation/assumptions because  the present analysis is limited (at most) to $T\in (0, T_{\rm f})$, but $T_{\rm f} =T_{1}= T_{\infty} < T_{\rm BD}$ when $\alpha = \alpha^{\bullet}$ is taken.

We conclude this section by recording the following limits, which correspond to the case 
$\rho_{\rm a}(x) \equiv \rho_{0}$:   
\be
\lim_{\alpha \to 0} \lshad P_{T} \rshad(T) = -\lim_{\alpha \to 0}\lshad P_{X} \rshad(T)=\frac{\pi}{1-\epsilon\hat{\beta}\pi T} \quad \implies \quad \lim_{\alpha \to 0}T_{\infty}=\gamma/(\epsilon\beta \pi);
\en
these expression should be compared with their counterparts in Ref.~\cite{JSV05}, wherein a different non-dimensionalization scheme was employed.


\section{Numerical scheme construction}\label{sect:Numerical}


\subsection{Overview of KSS methods}

Before providing the details of the proposed numerical methods that will be applied to IBVP~\eqref{IBVP:Atmos_Isothermal} and IBVP~\eqref{IBVP:LWE_non-D}, we will describe the 
Krylov subspace spectral (KSS) methodology in the context of its application to a simpler problem.
We consider the IBVP
\begin{subequations} \label{eq:KSSex}
\begin{align}
& u_{tt} + Lu = 0, \quad 0 < x < 1, \quad t > 0, \label{eq:KSSexPDE} \\
& u(0,t) = 0, \quad u(1,t) = 0, \quad t > 0, \label{eq:KSSexBC} \\
& u(x,0) = f(x), \quad u_t(x,0) = g(x), \quad 0 < x < 1, \label{eq:KSSexIC}
\end{align}
\end{subequations}
where the differential operator $L$ is defined by
\begin{equation}\label{eq:KSSexL}
Lu = -\kappa_{1}u_{xx} + q(x)u.
\end{equation}
We discretize on a uniform $N$-point grid $\{ x_i \}_{i=1}^N$, where $x_i = i\Delta x$ and $\Delta x = 
1/(N+1)$.  This yields the system of ODEs
\begin{subequations} \label{eq:KSSexODEs}
\begin{align}
& {\bf u}_N''(t) = -L_N {\bf u}_{N}(t), \quad t > 0, \\
& {\bf u}_N(0) = {\bf f}_N, \quad {\bf u}_N'(0) = {\bf g}_N,
\end{align}
\end{subequations}
where $L_N$ is a $N\times N$ matrix obtained via spatial discretization of $L$, say, by finite
differences, and ${\bf f}_N$ and ${\bf g}_N$ are discretizations of $f(x)$ and $g(x)$, respectively.  Moreover, in this section a prime indicates the derivative of a function of a single (real) variable. This  system of ODEs has the exact solution
\begin{subequations} \label{eq:KSSexsoln}
\begin{align}
{\bf u}_N(t+\Delta t) &= F_{11}(L_N, \Delta t) {\bf u}_N(t) + F_{12}(L_N, \Delta t) {\bf u}_N'(t), \\
{\bf u}_N'(t+\Delta t) &= F_{21}(L_N, \Delta t) {\bf u}_N(t) + F_{22}(L_N, \Delta t) {\bf u}_N'(t), 
\end{align}
\end{subequations}
where the matrix functions $F_{ij}(L_N, \Delta t)$, $i,j=1,2$, are given by
\begin{equation}
F_{11}(L_N, \Delta t) = F_{22}(L_N, \Delta t) = \cos(L_N^{1/2}\Delta t), \quad
F_{12}(L_N, \Delta t) = L_N^{-1/2}\sin(L_N^{1/2}\Delta t), \quad F_{21}(L_N, \Delta t) = 
-L_N^{1/2}\sin(L_N^{1/2}\Delta t).
\end{equation}
However, for large values of $N$, it is not practical to compute these matrix functions directly.

Let ${\bf u}_N^n$ represent an approximate solution of Sys.~\eqref{eq:KSSexODEs} at time $t_n$. A KSS method computes ${\bf u}_N^n$ as follows:
\begin{subequations} \label{eq:KSSapsoln}
\begin{align}
{\bf u}_N^{n+1} &= \tilde{F}_{11}(L_N, \Delta t) {\bf u}_N^n + \tilde{F}_{12}(L_N, \Delta t) [{\bf u}_N']^n, \\
[{\bf u}_N']^{n+1} &= \tilde{F}_{21}(L_N, \Delta t) {\bf u}_N^n + \tilde{F}_{22}(L_N, \Delta t) [{\bf u}_N']^n. 
\end{align}
\end{subequations}
In this paper, we will focus exclusively on a KSS method for systems of this form, with the addition of
source terms, that are second-order accurate in time; higher-order methods are described in 
Ref.~\cite{waveblock} for wave propagation problems, and in Refs.~\cite{kssepi,block} for other types of PDEs.

Let $S_N$ be the matrix of the $N$-point discrete sine transform. For $i,j=1,2$, $\tilde{F}_{ij}(L_N, \Delta t)$ is an approximation of $F_{ij}(L_N, \Delta t)$ that has the form
\begin{equation} \label{eq:FijKSS}
\tilde{F}_{ij}(L_N, \Delta t) = S_N^{-1} \left[  B^{ij}(\Delta t) S_N + M^{ij}(\Delta t) S_N L_N\right],
\end{equation}
where $B^{ij}(\Delta t)$ and $M^{ij}(\Delta t)$ are diagonal matrices.
For $k=1,2,\ldots,N$, the $k$th diagonal entries
$[B^{ij}(\Delta t)]_{kk}$ and $[M^{ij}(\Delta t)]_{kk}$ are the $y$-intercept and slope, respectively, of the
function $\mathcal{P}_{ij,k}(\lambda,\Delta t)$, that, for fixed $\Delta t$, is a linear function that interpolates $F_{ij}(\lambda,\Delta t)$ at frequency-dependent interpolation points 
$\lambda_{1,k},\lambda_{2,k}$.  


It remains to determine the interpolation points $\lambda_{1,k},\lambda_{2,k}$ for each
wave number $k$.  In this paper, the matrix $L_N$ is a finite difference discretization of the differential
operator $L$ in which the second derivative is approximated by a centered difference.  Therefore,
as in Ref.~\cite{resterms}, we choose 
\begin{equation} \label{eq:kssnodes}
\lambda_{1,k} = 0, \qquad \lambda_{2,k} = \kappa_{1} (2-2\cos(\pi k\Delta x))/(\Delta x)^2 + \avg(q),
\end{equation}
where $\avg(q)$ denotes the average value of the coefficient $q(x)$ on $(0,1)$.
That is, $\lambda_{2,k}$ is an approximate eigenvalue of $L_N$.


This approach to solving IBVP~\eqref{eq:KSSex} is second-order accurate in time and has
spectral accuracy in space; see Ref.~\cite{resterms} for a convergence analysis, including a proof that
the method is 
unconditionally stable for a similar problem,
in which periodic boundary conditions were imposed.  This kind of stability normally occurs only in
time-stepping methods that are either implicit, or involve some kind of iteration such as Krylov projection~\cite{HL97}.  

Krylov subspace methods for approximating matrix function-vector products of
the form $\varphi(A){\bf b}$ for a matrix $A$ and vector ${\bf b}$, such as those 
described in Ref.~\cite{HL97} in the case of the matrix exponential, 
generate a Krylov subspace of sufficient dimension to ensure
the desired accuracy of the approximation.  When such methods are used to solve stiff
systems of ODEs derived from PDEs, the required Krylov subspace dimension can grow
substantially when increasing the time step, or, depending on the PDE, the number of grid
points in the spatial discretization \cite{kssepi}.  

By contrast, in a KSS method, the
Krylov subspace dimension is independent of the time step or spatial grid size; rather, it is
determined solely by the order of temporal accuracy.  For example, in the KSS method
defined by Eq.~\eqref{eq:FijKSS}, for each approximation of a matrix function-vector product, 
only one matrix-vector multiplication (excluding multiplication by diagonal matrices, which of course is implemented as component-wise multiplication of vectors), at most two Fourier transforms, and one inverse Fourier transform
are required.  As such, KSS methods
realize stability
and scalability~\cite{kssepi} superior to that of similar time-stepping methods, through the use of
frequency-dependent approximations of the solution operator of the PDE.  In the remainder of this
section, we discuss how a second-order KSS method can be applied to the 
problems presented in this paper.

\subsection{Shock case} \label{sect:shock_scheme}

To solve  hIBVP (\ref{IBVP:Atmos_Isothermal}) numerically, we first modify the problem so that it is defined on a bounded domain $(0,\ell)$. That is, we  consider
\begin{subequations} \label{IBVP:eqorigshock} \begin{align}
&w_{tt}-c_0^2 w_{zz}+\gamma g w_z+\hat{\mu} w_t=0, \quad 0 < z < \ell, \quad t > 0, 
\label{eq:EoM_eqorigshock}\\
&w(z,0) = w_t(z,0) = 0 ,\\
&w(0,t)={\cal H}(t)\cos(\omega t),
\end{align}
\end{subequations}
where, for the moment, we have set $W_0=1$. Since the boundary condition causes a shock wave that is entering from the left, the solution will be zero to the right of the wave front. Therefore, we impose the right-end boundary condition $w(\ell,t)=0$. Due to the Dirichlet boundary conditions, we use a Fourier sine series. Thus, we must transform the problem to eliminate first derivative terms and homogenize the the left boundary condition. 

\subsubsection{Transforming the PDE} \label{sect:shock_transformpde}

We first simplify  Eq.~\eqref{eq:EoM_eqorigshock}, which we note appears (in 
IBVP~\eqref{IBVP:eqorigshock}) unchanged from its role as  the EoM  in  hIBVP~\eqref{eq:Atmos_Ray_EoM}, by introducing a differential operator $L$, defined by
$$
Lw=-a_2w_{zz}+a_1w_z,
$$
where $a_2=c_0^2$ and $a_1=\gamma g$; in terms of $L$, then, Eq.~\eqref{eq:EoM_eqorigshock} becomes
\begin{equation} \label{eqorigPDE} w_{tt} = -Lw-\hat{\mu}w_t.\end{equation}
Proceeding as suggested in Ref.~\cite{517txtbook}, by defining 
$$w(z,t)=\psi(z,t)\tilde{w}(z,t), \quad \psi(z)=\exp\left[\left(\frac{a_1}{2a_2}\right) z\right],$$
we obtain the equivalent PDE
$$\tilde{w}_{tt}=-\tilde{L}\tilde{w}-\hat{\mu}\tilde{w}_t,$$
where $$\tilde{L}\tilde{w}=-a_2 \tilde{w}_{zz}+ a_0 \tilde{w}, \quad a_0=\frac{a_1^2}{4a_2}.
$$

\subsubsection{Homogenizing the boundary conditions} \label{sect:shock_homobc}

We now transform the hIBVP into an equivalent problem with homogeneous boundary conditions,
\begin{subequations}  \label{IBVP_hom} \begin{align} 
&u_{tt}=-\tilde{L}u-\hat{\mu}u_t-G(z,t), \label{PDE_hom} \\
&u(0, t) = 0, \, u(\ell,t)=0, \\
&u(z,0) = -F(z,0), \, u_t(z,0) = -F_t(z,0),
\end{align} \end{subequations}
where 
$$G(z,t) = F_{tt} + a_2 F_{zz} + a_0 F + \hat{\mu}F_{t}$$ 
is the source term introduced by transforming the left boundary condition.  We use
an approach to homogenization of boundary conditions
described in Ref.~\cite{517txtbook}, except that we also require that the source term $G(z,t)$ satisfies 
\begin{subequations} \begin{align}
& G(0,t)=0, \label{condG0}\\
& G(z,t)=0, \quad z\geq \tilde{\ell}, \label{condGell}
\end{align} \end{subequations} 
for some $\tilde{\ell}<\ell$.  We choose $\tilde{\ell}=80\Delta z$ for the implementation, where $\Delta z$ is the spatial grid mesh, but numerical experimentation indicates that the performance of the numerical method is not sensitive to the choice of $\tilde{\ell}$.  A similar approach to homogenization was used in Ref.~\cite{nmpde}, except that in this case, we require that $G(z,t)$ vanishes on most of the spatial domain, instead of only at the right endpoint.  This mitigates numerical artifacts, caused by truncation error, in the solution to the right of the wavefront.

The solution of Eq.~\eqref{IBVP_hom} is then written as 
$$u(z,t)=\tilde{w}(z,t)-F(z,t),$$
where 
$$
F(z,t) = \begin{cases}
f_0(t) + f_1(t) z + f_2(t) z^2 + f_3(t) z^3 + f_4(t) z^4 & 0 \leq z \leq \tilde{l} \\
0 & \tilde{l} < z \leq \ell.
\end{cases}
$$
satisfies not only the
original boundary conditions, but also belongs to $C^2([0,\ell])$.  That is, we require
\begin{subequations} \begin{align}
& F(0,t)=\cos(\omega t), \\
& F(\tilde{\ell},t)=0, \label{condFell}\\
& F_z(\tilde{\ell},t)=0, \label{condF_zell} \\
& F_{zz}(\tilde{\ell},t)=0. \label{condF_zzell} 
\end{align} \end{subequations} 
The smoothness requirements on $F$, an additional modification of the homogenization scheme
from Ref.~\cite{nmpde}, are imposed to increase the decay rate of Fourier coefficients of the solution
and thus reduce high-frequency oscillations.

From condition~(\ref{condG0}), we obtain 
\begin{eqnarray*}
 f_2(t) &=& \frac{1}{2a_2}\left(f_0''(t) + a_0 f_0(t) +\hat{\mu}f_0'(t)\right) \\
 &=& -\frac{\omega^2}{2a_2}\cos(\omega t) + \frac{a_0}{2a_2} \cos(\omega t) - \frac{\hat{\mu} \omega}{2a_2}\sin(\omega t). 
\end{eqnarray*} 
We then use the conditions  (\ref{condGell}), (\ref{condFell}), and (\ref{condF_zell}) in a $3\times 3$ system of linear equations to obtain $f_1(t)$, $f_3(t)$, and $f_4(t)$.  
Condition~(\ref{condF_zzell}) is 
also satisfied as a result.

\subsubsection{Approximate solution operator} \label{sect:shock_solnop}

Now, we express the solution at time $t_{n+1}$ in terms of the solution at time $t_n=n\Delta t$.
We will work with a PDE of the form
\begin{equation} \label{eqPDEst}
u_{tt} = -\tilde{L}u + b(z,t_n),
\end{equation} 
where $\tilde{L}$ is a positive definite second-order spatial differential operator, 
and $b(z,t_n)=-\hat{\mu}u_t(z,t_n)-G(z,t_n)$ is treated as a source term, as a first-order approximation of Eq.~\eqref{PDE_hom}. We then rewrite Eq.~\eqref{eqPDEst} as a first-order system
\begin{eqnarray} \label{eqsyst}
\left[ \begin{array}{c}
u \\ u_t \end{array} \right]_t & = & 
J \left[ \begin{array}{c} u \\ u_t \end{array} \right] + 
\left[ \begin{array}{c} 0 \\ b \end{array} \right],\quad J=\left[ \begin{array}{cc}
0 & \mathbb{I} \\ -\tilde{L} & 0 \end{array} \right],
\end{eqnarray}
where we use $\mathbb{I}$ to denote the identity operator.
If we define
$$
{\bf r} = \left[ \begin{array}{c} u \\ u_t \end{array} \right], \quad
{\bf c}_n = \left[ \begin{array}{c} 0 \\ b(z,t_n) \end{array} \right],
$$
then, as in Ref.~\cite{ostermann22}, we can express the solution of Sys.~\eqref{eqsyst} as
\begin{equation} \label{solnoperator}
{\bf r}(t_{n+1}) = \phi_0(J\Delta t) {\bf r}(t_n) + \Delta t\phi_1(J\Delta t) {\bf c}_n, 
\end{equation}
where 
\begin{equation} \label{eq:phi0}
\phi_0(Jt) = 
\left[ \begin{array}{cc}
\cos(\tilde{L}^{1/2} t) & \tilde{L}^{-1/2} \sin(\tilde{L}^{1/2} t) \\
-\tilde{L}^{1/2} \sin(\tilde{L}^{1/2} t) & \cos(\tilde{L}^{1/2} t)
\end{array} \right]
\end{equation}
and
\begin{equation}
t\phi_1(Jt) 
 =  \left[ \begin{array}{cc}
\tilde{L}^{-1/2} \sin(\tilde{L}^{1/2} t) &\tilde{L}^{-1}-\tilde{L}^{-1} \cos(\tilde{L}^{1/2} t) \\ 
\cos(\tilde{L}^{1/2} t) - \mathbb{I} & \tilde{L}^{-1/2} \sin(\tilde{L}^{1/2} t)
\end{array} \right]. \label{eq:phi1}
\end{equation}
What remains is to apply spatial discretization to Eq.~\eqref{solnoperator}, and then
compute the various matrix function-vector products indicated by Eqs.~\eqref{eq:phi0} and~\eqref{eq:phi1}.

\subsubsection{Numerical solution via KSS} \label{sect:shock_KSS}

The computed solution (Eq.~\eqref{solnoperator}) of the PDE involves approximating functions of a matrix $\tilde{L}_N$,
that is identified with a differential operator $\tilde{L}$ that $\tilde{L}_N$ discretizes through centered differencing.  The operator $\tilde{L}$ in this case will be
the result of the transformation described previously,
$$
\tilde{L}u = -a_2 u_{zz} + a_0 u.
$$ 
To approximate $f(\tilde{L}_N){\bf u}_{N}$ for a given function $f$ and vector ${\bf u}_{N}$ using a second-order KSS approximation, as in 
Sys.~\eqref{eq:KSSapsoln}, 
we first select frequency-dependent interpolation points as in
Eq.~\eqref{eq:kssnodes}, 
\begin{equation} \label{eq:shocknodes}
\lambda_{1,k} = 0, \qquad \lambda_{2,k} = a_2 (2-2\cos(\omega_k\Delta z))/(\Delta z)^2 + a_{0},
\end{equation}
where $\omega_k=\pi k/\ell$, for $k=1,2,\ldots,N$, and $\lambda_{2,k}$ is an eigenvalue of $\tilde{L}_N$.
Based on Eq.~\eqref{eq:FijKSS}, the approximation is then
\begin{equation} \label{eq:KSSapprox}
f(\tilde{L}_N) {\bf u}_{N}\, \approx \,S_N^{-1} [ B_k   S_N {\bf u}_{N} + M_k  S_N \tilde{L}_N{\bf u}_{N}],
\end{equation}
where $B_k$ and $M_k$ are diagonal matrices, the diagonal entries of which are 
the $y$-intercepts and slopes, respectively, of the linear approximation of $f$ for each frequency $\omega_k$ with interpolation points from Eq.~\eqref{eq:shocknodes}.  The operators $S_N$ and $S_N^{-1}$ represent the $N$-point discrete sine transform and its inverse, respectively.

It is worth noting that because the operator $\tilde{L}$ has constant coefficients, 
Eq.~\eqref{eq:KSSapprox} can be simplified to 
\begin{equation} \label{eq:FSapprox}
f(\tilde{L}_N){\bf u}_{N}\approx S_N^{-1} f(\Lambda_2) S_N {\bf u}_{N},
\end{equation}
where $\Lambda_2$ is a diagonal matrix with diagonal entries $\lambda_{2,k}$, $k=1,2,\ldots,N$.
In this case, the KSS method reduces to a Fourier spectral method.
However, in numerical experiments, it was found that this simplified formula led to high-frequency
oscillations in the solution, known as Gibbs' phenomenon, brought on by substantial accumulation of round-off error, whereas using Eq.~\eqref{eq:KSSapprox} as written did not.  As will be illustrated
later in this section, even with enough regularization to keep the solution bounded, using a Fourier 
spectral method, even in this constant-coefficient case, is not a viable approach.

\subsection{Lighthill--Westervelt equation} \label{sect:LWE_scheme}

The process of transforming IBVP~\eqref{IBVP:LWE_non-D} is similar to that of the shock case, 
except that additional steps are required due to the spatial variation of the coefficients.  For
convenience, we recast  IBVP~\eqref{IBVP:LWE_non-D} as follows:
\begin{subequations}\label{IBVP:LWE_short}
\begin{align}
\label{eq:EoM_P_short}
&P_{TT} + {\cal L}P - \epsilon\hat{\beta} \partial_{TT}(P^2) = 0,  \quad  (X,T)\in (0, 1)\times (0, T_{\rm f}),\\
&P(0,T) = \mathcal{H} (T)\sin(\pi T),\quad P(1,T) =0, \qquad  T\in (0, T_{\rm f}), \\
&P(X,0) = 0, \quad  P_{T}(X,0) = 0, \qquad  X\in (0,1),
\end{align}
\end{subequations}
where the differential operator ${\cal L}$ is defined by
$${\cal L}P = -V_{\rm a}^2(X) P_{XX} + \alpha V_{\rm a}^2(X) P_X,$$
with 
$$V_{\rm a}^2(X) = \exp(2\mathfrak{K} X), $$
and where we let  $\mathfrak{K}  =-\alpha/2$ for convenience.

\subsubsection{Transforming the PDE}

As in \cite{nmpde}, to homogenize the coefficient of $P_{XX}$, we define a change of independent variable 
$$
Y=\phi(X) = \frac{1-e^{-\mathfrak{K} X}}{1-e^{-\mathfrak{K} }}.
$$
Then, we obtain the transformed differential operator $\tilde{\cal L}$, defined by
$$
\tilde{\cal L} U(Y) = -\tilde{a}_2 U_{YY} + \tilde{a}_1(Y) U_Y,
$$
where
$$
\tilde{a}_2 = C^2, \qquad \tilde{a}_1(Y) = (\mathfrak{K} +\alpha)C \exp[\mathfrak{K} \phi^{-1}(Y)],
$$
with
$$C = \left[ \int_{0}^{1} \exp(-\mathfrak{K}X)\,\rd X \right]^{-1} = \frac{\mathfrak{K}}{1 - \exp(-\mathfrak{K})},$$
being the homogenized wave speed.

Next, as in Section~\ref{sect:shock_transformpde}, we eliminate the coefficient of $U_Y$,
which yields the final form of our spatial differential operator,
\begin{equation} \label{eq:barLdef}
\bar{\cal L} U = -\bar{a}_2 U_{YY} + \bar{a}_0(Y) U,
\end{equation}
where
$$
\bar{a}_2 = C^2, \qquad \bar{a}_{0}(Y) =  \frac{\tilde{a}_1(Y)^2}{4\tilde{a}_2} 
- \frac{\tilde{a}_1'(Y)}{2}.
$$
The transformed problem is now
\begin{subequations}\label{IBVP:LWE_trans}
\begin{align}
\label{eq:EoM_P_trans}
&\bar{P}_{TT} +  \bar{\cal L}\bar{P} - 2\epsilon\hat{\beta}\partial_{T}( \bar{P}\bar{P}_T) = 0,  \quad  (Y,T)\in (0, 1)\times (0, T_{\rm f}),\\
&{\bar P}(0,T) = \mathcal{H} (T)\sin(\pi T),\quad \bar{P}(1,T) =0, \qquad  T\in (0, T_{\rm f}), \\
&\bar{P}(Y,0) = 0, \quad  \bar{P}_T(Y,0)  = 0, \qquad  Y\in (0,1),
\end{align}
\end{subequations}
where
$$\bar{P}(Y,T) = \psi^{-1}(Y) P(\phi^{-1}(Y),T), \quad 
\psi(Y) = \exp\left[ \int \frac{\tilde{a}_1(Y)}{2\tilde{a}_2}\,\rd Y \right].
$$

\subsubsection{Homogenizing the boundary conditions}

Next, as in Section~\ref{sect:shock_homobc}, we homogenize the boundary conditions.  If we let
$$
\bar{P} = U+F,
$$
then IBVP~\eqref{IBVP:LWE_trans} becomes
\begin{subequations}\label{IBVP:LWE_BC}
\begin{align}
\label{eq:EoM_P_BC}
&{U}_{TT} + \bar{\cal L}{U} - 2\epsilon\hat{\beta} \partial_{T}\{ UU_T + U_T F + UF_T \} = G,  \quad  (Y,T)\in (0, 1)\times (0, T_{\rm f}),\\
&U(0,T) = 0,\quad \quad U(1,T) = 0, \qquad  T\in (0, T_{\rm f}), \\
&U(Y,0) = -F(Y,0), \quad  {U}_T(Y,0) = -F_T(Y,0), \qquad  Y\in (0,1),
\end{align}
\end{subequations}
where
$$
G=-{F}_{TT} - {\cal L}F + 2\epsilon\hat{\beta} \partial_{T}\{F F_T \}.
$$
The functions $F(Y,T)$ and $G(Y,T)$ satisfy the following conditions for $T>0$:
\begin{itemize}
\item $F(0,T) = {\cal H}(T)\sin(\pi T)$, $G(0,T) = 0,$
\item For $Y\geq \ell$, $F(Y,T) = 0$, $F_Y(Y,T) = 0$, $G(Y,T) = 0,$
where $\ell$ is chosen to be $1/2$ in our implementation.
\end{itemize}
The computation of $F$ and $G$ proceeds as described in Section~\ref{sect:shock_homobc}.


\subsubsection{Numerical solution via KSS} \label{subsubsec:LWEnumKSS}

A second-order KSS method is then applied to IBVP~\eqref{IBVP:LWE_BC}, almost exactly as described 
in Sections~\ref{sect:shock_solnop} and~\ref{sect:shock_KSS}.  That is, we use the KSS method
to approximate 
\begin{equation} \label{eq:LWEsolnoperator}
{\bf r}(T_{n+1}) = \phi_0(J\Delta T) {\bf r}(T_n) + \Delta T\phi_1(J\Delta T) {\bf c}_n, 
\end{equation}
where $T_n=n\Delta T$ and 
$$
J = \left[ \begin{array}{cc}
0 & \mathbb{I} \\
-\bar{\cal L} & 0 
\end{array} \right], \qquad
{\bf r} = \left[ \begin{array}{c} U \\ U_T \end{array} \right], \qquad
{\bf c}_n = \left[ \begin{array}{c} 0 \\ b(Y,T_n) \end{array} \right],
$$
to solve the PDE
\begin{equation} \label{eq:LWEPDEst}
U_{TT} = -\bar{\cal L}U + b(Y,T_n),
\end{equation} 
where
\begin{equation} \label{eq:LWEbdef}
b(Y,T_n) = G(Y,T_n) + 2\epsilon\hat{\beta} Q(Y,T_n), \qquad
Q = \partial_T\{ UU_T + U_T F + UF_T \}.
\end{equation}
One substantial difference from the shock case is that
we also need to discretize $U_{TT}$ to compute $Q$.  At the $n$th time step, this is accomplished by 
$$
{\partial}_{TT} {\bf U}_N^n \approx S_N^{-1} \left[ \frac{\partial_T \hat{\bf U}_N^n - \partial_T \hat{\bf U}_N^{n-1}}{\Delta T} \right], \quad \hat{\bf U}_N^n = S_N {\bf U}_N^n.
$$
Prior to taking the inverse sine transform, the above difference quotient is regularized to eliminate
Gibbs' phenomenon, as is done for the time derivative of the solution.  This is discussed further
in Section~\ref{sect:Numerical_LWE}.

The same KSS method described in Section~\ref{sect:shock_KSS} is used to approximate
Eq.~\eqref{eq:LWEsolnoperator}.  This calls for approximating matrix functions of $\bar{\cal L}_N$, a
$N\times N$ matrix that discretizes the operator $\bar{\cal L}$ via centered differencing.  For this
task, the formula given in Eq.~\eqref{eq:KSSapprox} is again employed, with interpolation points
\be
\lambda_{1,k} = 0, \quad \lambda_{2,k} = a_{2} (2-2\cos(\pi k\Delta Y))/(\Delta Y)^2 + \avg(\bar{a}_0), \qquad k =1,2,\ldots,N,
\en
where $\avg(\bar{a}_0)$ denotes the average value of $\bar{a}_{0}(Y)$ on $(0,1)$.

\subsection{Numerical results for the shock wave case}\label{sect:Numerical_Shock}

\begin{figure}[ht] 
\includegraphics[width=2.2in]{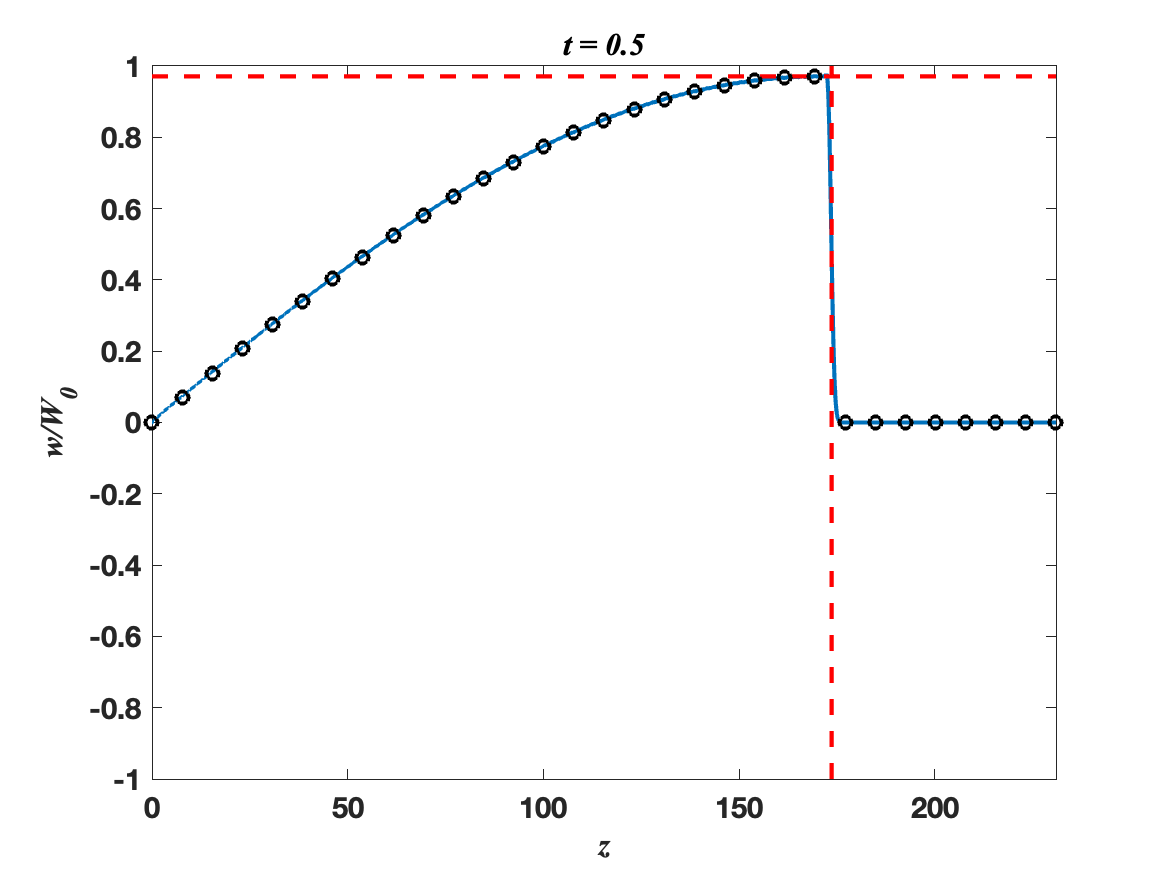}
\includegraphics[width=2.2in]{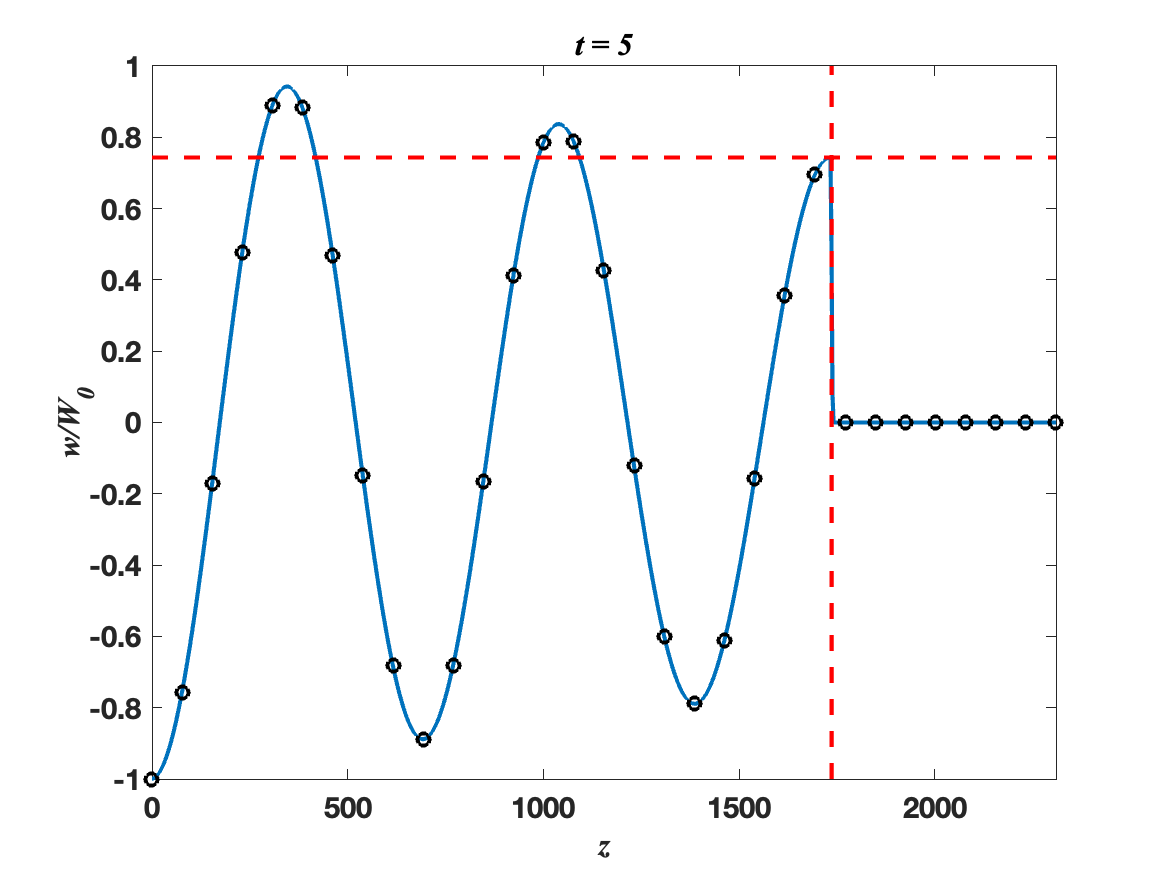}
\includegraphics[width=2.2in]{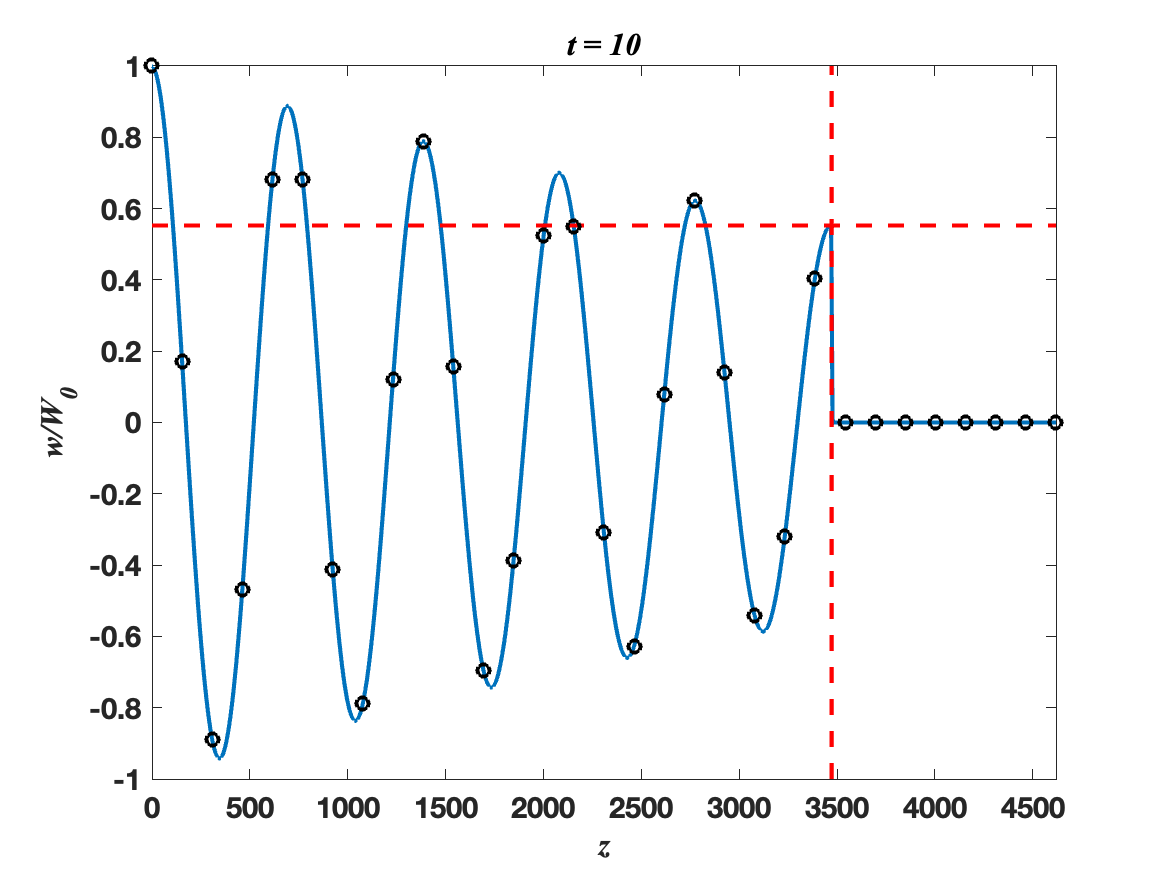}
\caption{Blue curves:  Numerical solution of  IBVP~\eqref{IBVP:eqorigshock} computed using the KSS method
described in Section~\ref{sect:shock_scheme} with $N=262,144$ grid points and CFL number 16.3841. The black circles represent the analytical solution; they were generated by numerically evaluating  Eq.~\eqref{eq:Ray_t-domain-sol} using \textsc{Mathematica's}
 (ver.~11.2)   \texttt{NIntegrate[\,$\cdot$\,]} command. The horizontal and vertical red-dashed lines represent the theoretical shock amplitude values (see Eq.~\eqref{eq:shock-Amp_Isothermal})  and wave-front locations (see Eq.~\eqref{eq:shock-Loc_Isothermal}), respectively.}\label{fig:Shock}
\end{figure}

The sequence shown in Fig.~\ref{fig:Shock} depicts particular instances in the evolution of the $W$ vs.\ $z$ solution profile in the case of IBVP~\eqref{IBVP:Atmos_Isothermal}, where we have set $W :=w/W_{0}$; it is based on the following parameter values, which, if not dimensionless,  are stated in SI units:
\begin{multline}
\vartheta_{0} = 300, \quad   \gamma = 1.40, \quad c_{0}=347.26, \quad  \omega = 2\pi,  \quad g=9.81, \quad \hat{\mu}= 2\hat{\mu}_{\rm c} \approx 0.0791,
\end{multline}
where the values taken for $\gamma$, $g$, and $c_{0}$ yield $\hat{\mu}_{\rm c} \approx 0.0395$.
We set the right endpoint $\ell$ of the spatial domain in Eq.~\eqref{IBVP:eqorigshock} to be $40c_0$.

Fig.~\ref{fig:Shock} shows the solution obtained using $N=262,144$ grid points and a CFL number of $\frac{c_0\Delta t}{\Delta z}\approx 16.3841$ at times $t=0.5, 5, 10$.  To reduce oscillations due to the Gibbs phenomenon that occurs due to the discontinuity, the Fourier sine coefficients $W$ and $W_t$ are multiplied by Lanczos sigma factors~\cite{hamming}. The sigma factors are raised to appropriate powers $\mathfrak{p}$ and $\mathfrak{p}^{\diamond}$ for $W$ and $W_t$, respectively, to avoid the weak instability that would occur without any smoothing~\cite{GHT94}. Less regularization is needed over time. In the implementation $\mathfrak{p}=(1-t/t_{\rm end})^2$ and $\mathfrak{p}^{\diamond}=3(1-t/t_{\rm end})^2$, where $t_{\rm end}=11.0$ is the final time. The numerical solution very closely matches the analytical solution given in Eq.~\eqref{eq:Ray_t-domain-sol}. The computed solution also has the correct wave-front location (see Eq.~\eqref{eq:shock-Loc_Isothermal}) and the correct shock amplitude (see Eq.~\eqref{eq:shock-Amp_Isothermal}) at the wave front, as indicated by the vertical and horizontal 
red-dashed lines, respectively. Since the numerical solution is also vertical at the wave-front, it correctly captures the jump discontinuity, as desired.  

\begin{figure}[ht] 
\includegraphics[width=2.2in]{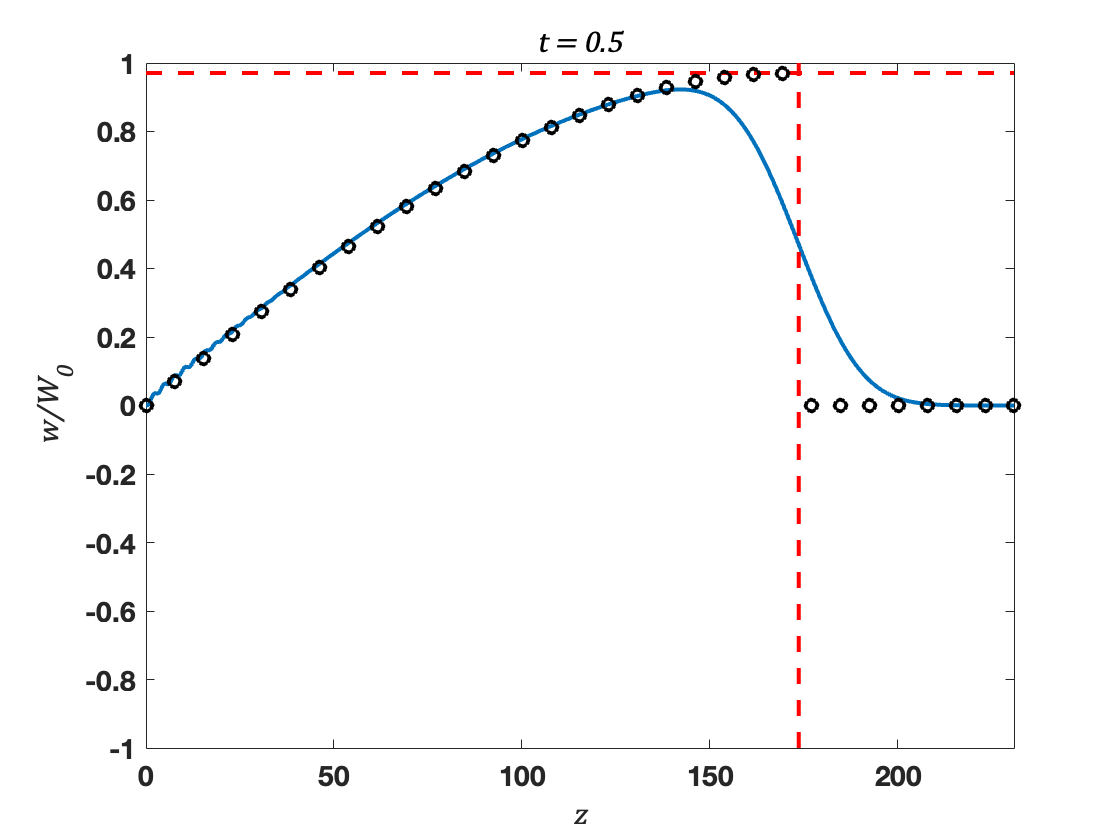}
\includegraphics[width=2.2in]{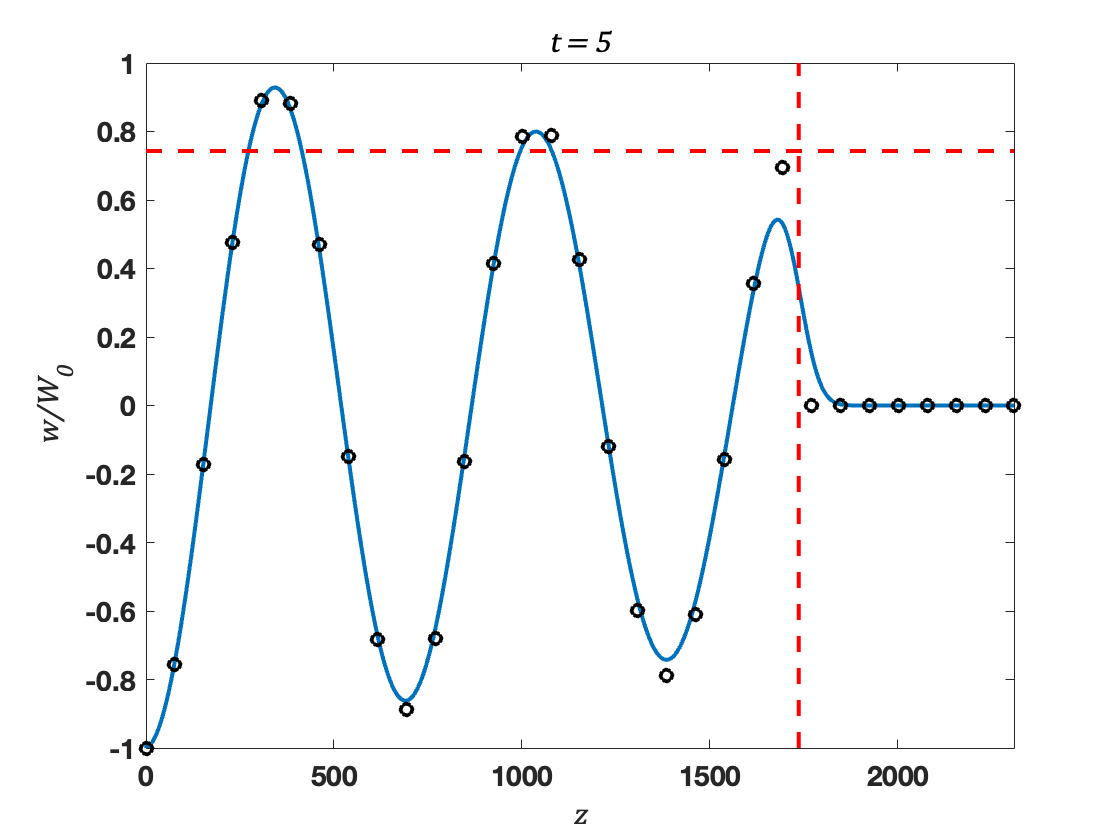}
\includegraphics[width=2.2in]{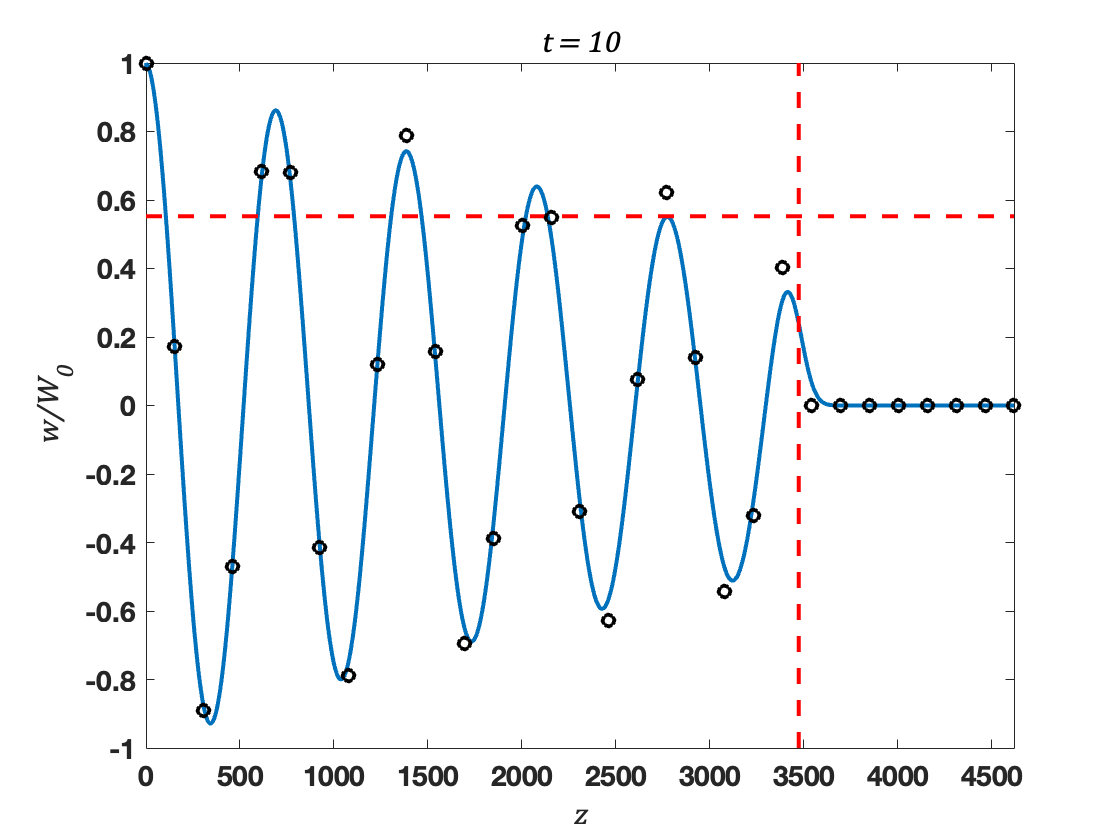}
\caption{Blue curves:  Numerical solution of  IBVP~\eqref{IBVP:eqorigshock} computed using the Fourier spectral method
described in Section~\ref{sect:shock_scheme} with $N=262,144$ grid points and CFL number 16.3841,
with Eq.~\eqref{eq:KSSapprox} replaced by Eq.~\eqref{eq:FSapprox}. The black circles represent the analytical solution; they were generated by numerically evaluating  Eq.~\eqref{eq:Ray_t-domain-sol} using \textsc{Mathematica's}
 (ver.~11.2)   \texttt{NIntegrate[\,$\cdot$\,]} command. The horizontal and vertical red-dashed lines represent the theoretical shock amplitude values (see Eq.~\eqref{eq:shock-Amp_Isothermal})  and wave-front locations (see Eq.~\eqref{eq:shock-Loc_Isothermal}), respectively.}\label{fig:ShockF}
\end{figure}

We then repeat the solution process, with two modifications. First, the formula for
matrix function-vector products from Eq.~\eqref{eq:KSSapprox} is simplified to 
the formula from Eq.~\eqref{eq:FSapprox}.  That is, a Fourier spectral method is used to compute the solution, instead of a KSS method.  Second, the regularization
applied to the KSS method is again used, but with the sigma factors raised to powers $\mathfrak{p}=
\mathfrak{p}^{\diamond} = 192$, as anything less results in instability.  We see
in Fig.~\ref{fig:ShockF} that the solution is excessively smoothed, so we conclude that this approach cannot produce an accurate solution, at least not without resorting to such effort to fine-tune the regularization that this approach is rendered impractical.

\subsection{Numerical results for the Lighthill--Westervelt equation case}\label{sect:Numerical_LWE}

%

\begin{figure}[ht]
\includegraphics[width=2.2in]{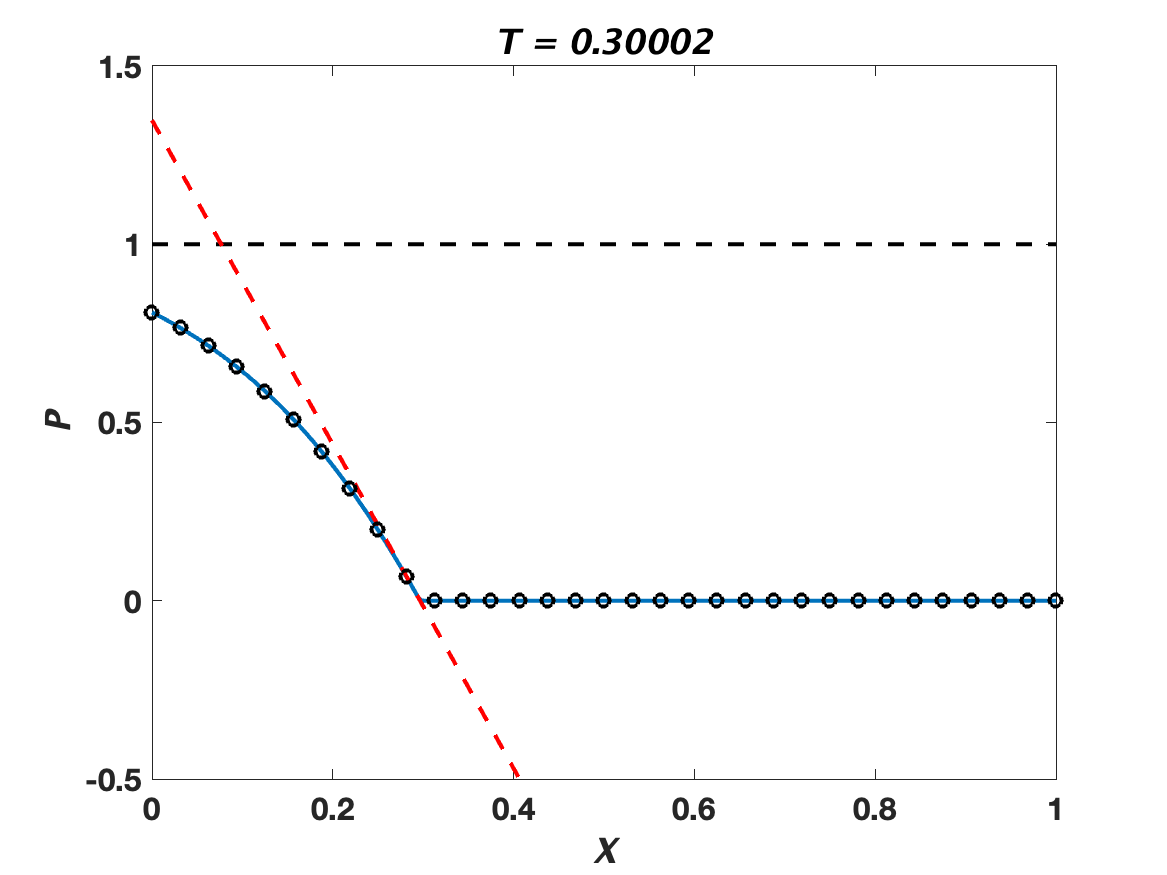}
\includegraphics[width=2.2in]{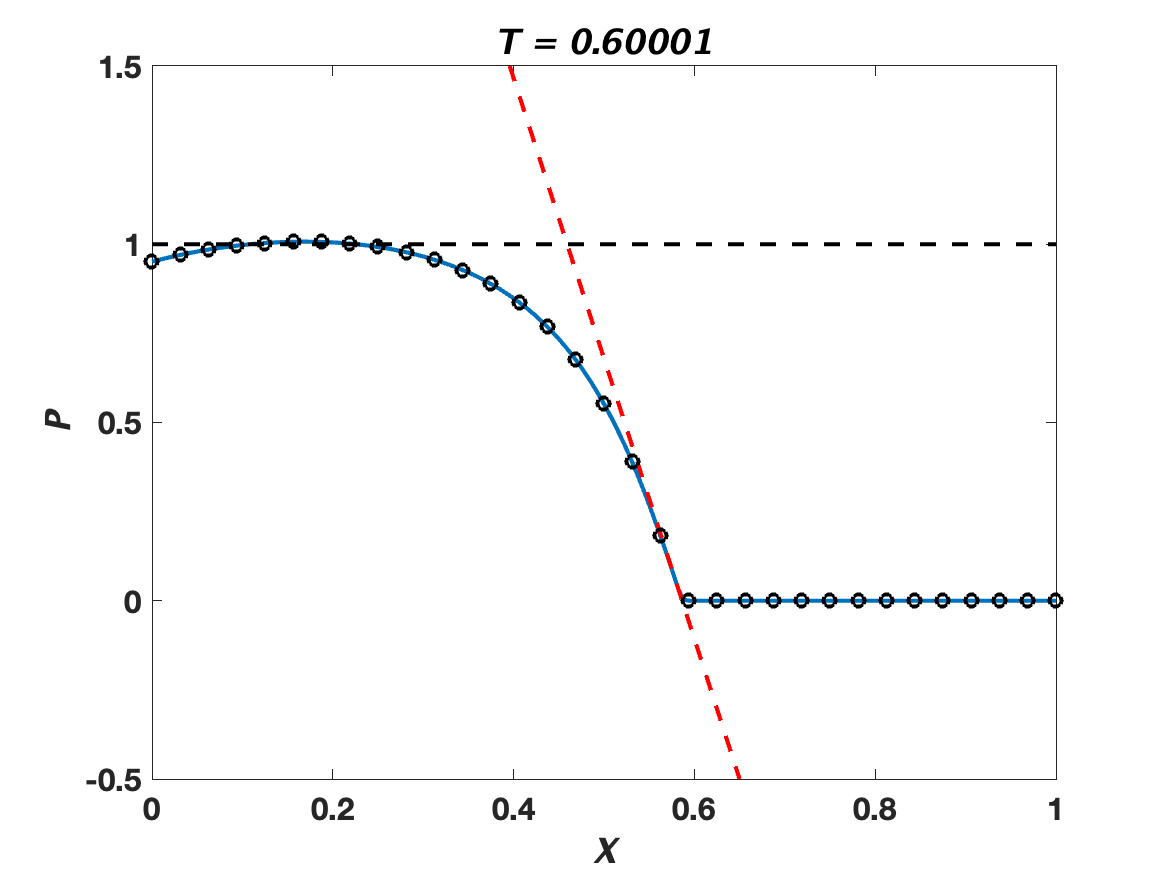}
\includegraphics[width=2.2in]{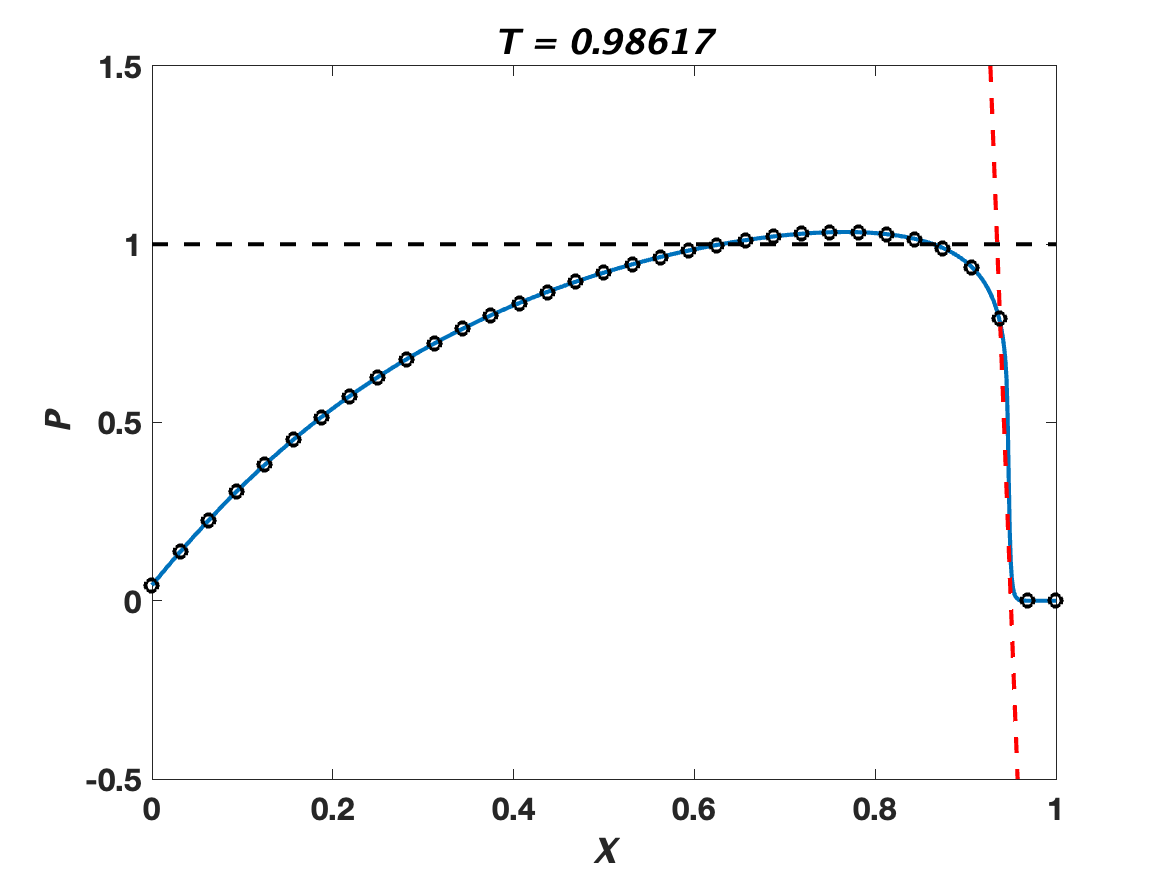}
\caption{Blue curves: Numerical solution of IBVP~\eqref{IBVP:LWE_short}  based on the parameters given in Eq.~\eqref{eq:LWE1parms},  computed using the KSS method described in Section~\ref{sect:LWE_scheme} with 
$N=262,144$ grid points and CFL number 10. The black circles represent the solution computed using FDS~\eqref{eq:FDS_P} with $M=4,096$  and CFL number $\bull \approx 0.5201$ (see Eq.~\eqref{eq:FDS_CFL}).  The red-dashed lines are plots of  $(X-\ram(T))\lshad P_{X} \rshad(T)$; they represent the tangent lines to  the solution profiles  at the wave-front $X=\ram(T)$.}
\label{figLWEp}
\end{figure}

The sequence shown in Fig.~\ref{figLWEp} depicts the solution of IBVP~\eqref{IBVP:LWE_short} based on 
the following \emph{dimensionless} parameter values:
\begin{equation}  \label{eq:LWE1parms}
\gamma = 1.40, \qquad \beta = 1.2, \qquad \epsilon = 0.35, 
\end{equation}
with $\alpha=\alpha^\bullet$, where $\alpha^\bullet$ is as defined in Eq.~\eqref{eq:alpha_bullet}.  This yields
$\alpha=\alpha^{\bullet}\approx 0.156451$ and $T_{\rm f}=T_{\infty} = T_{1} \approx 1.04015$.
The solution is plotted at $T\approx 0.3, 0.6, 0.98617$, where $X = \ram (0.98617)\approx 0.95$

The blue curves
were obtained via KSS using $N=262,144$ grid points and a CFL number of $\frac{\Delta T}{\Delta X}=10$.
Oscillations due to the Gibbs phenomenon are addressed as in the shock wave case, except that the regularization is applied
to $P_T$ and $P_{TT}$, whereas $P$ does not need to be regularized as it remains continuous.
Also, unlike the shock wave case, the Lanczos sigma factors are raised to a power that {\em increases}
with time, due to the `shocking-up' that occurs.  
More precisely, for the indicated grid size and
CFL number, the Lanczos sigma factors
for the Fourier sine transforms of $P_T$ and $P_{TT}$ are raised to powers 
$\varpi^{\diamond}$ and $\varpi^{\diamond\diamond}$, respectively,
where $\varpi^{\diamond} =1536(T/T_{\rm end})$ and $\varpi^{\diamond\diamond}=16384(T/T_{\rm end})^2$, where $T_{\rm end}=0.98617$.  
It can be seen that the slopes of the KSS solution profiles at the wave-fronts match the analytically-derived values (see Eq.~\eqref{eq:[PX]}).

\begin{figure}[ht]
\includegraphics[width=2.2in]{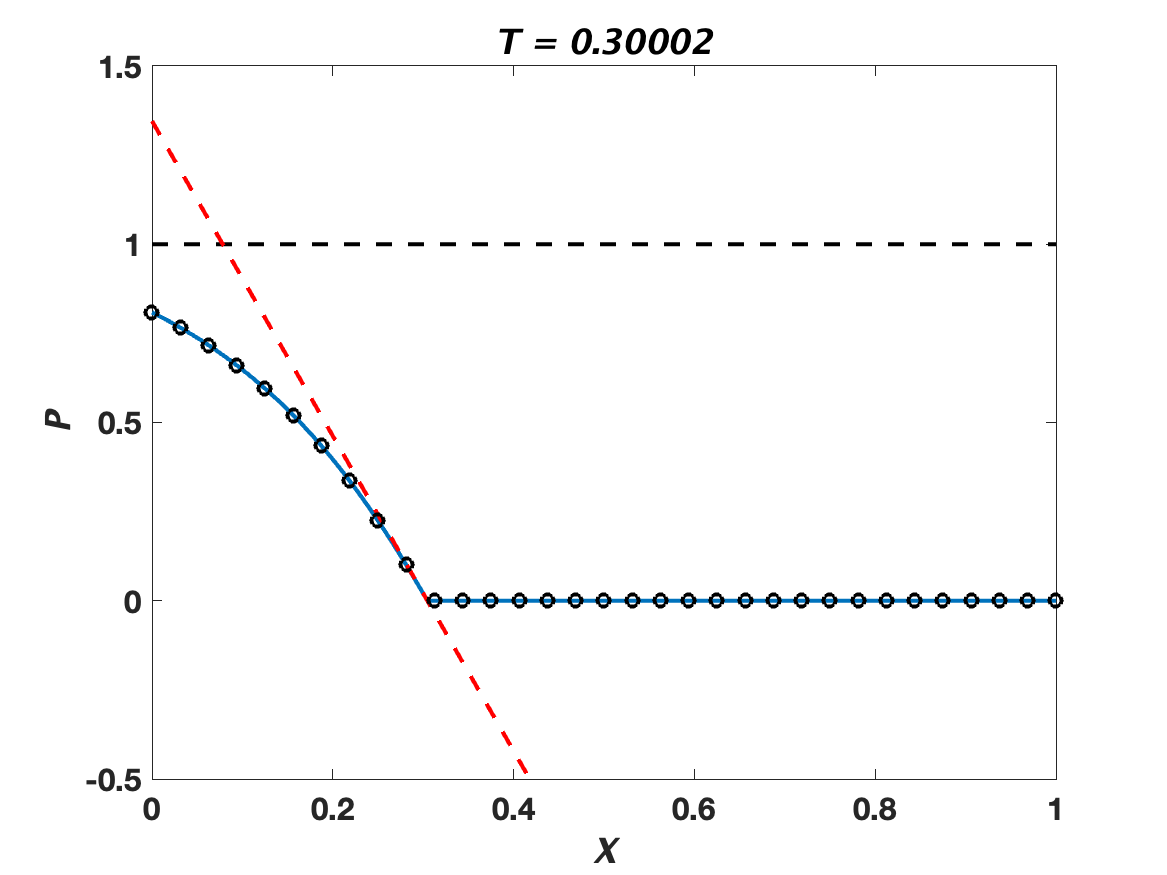}
\includegraphics[width=2.2in]{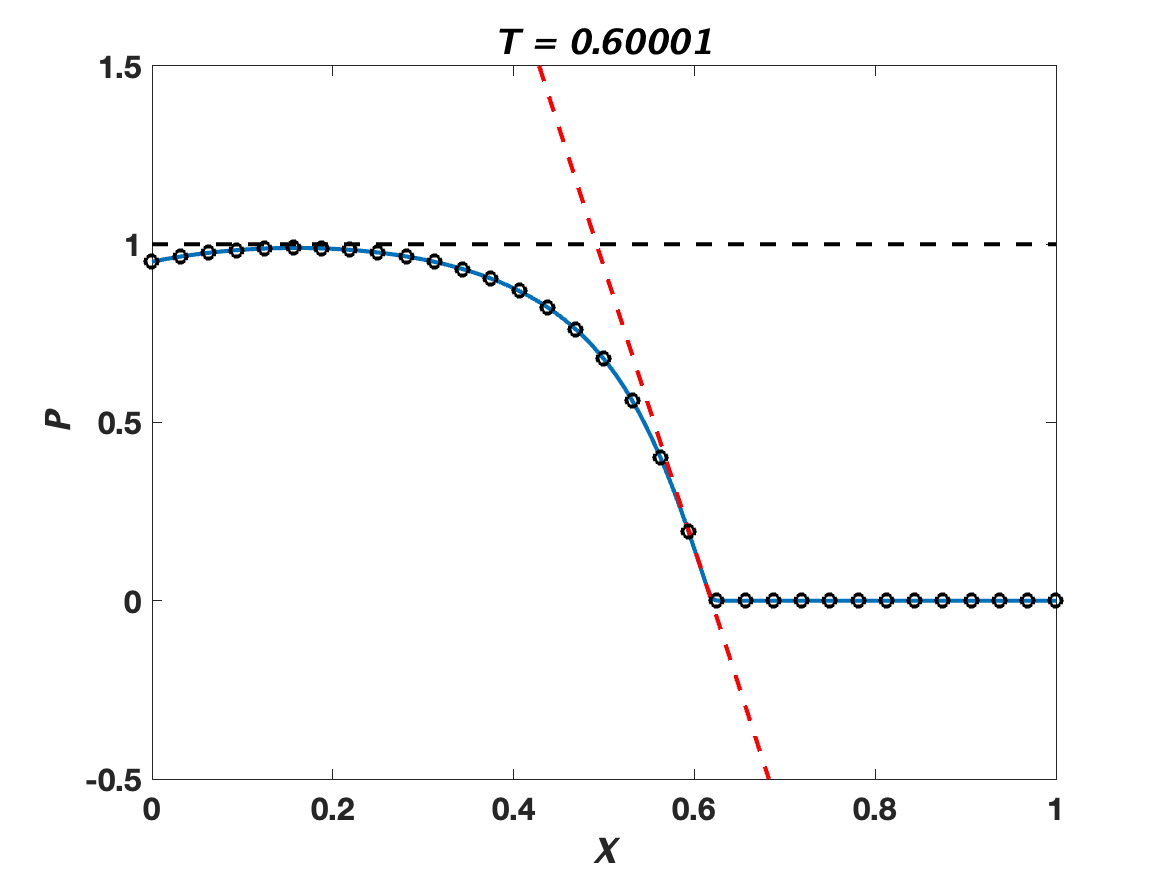}
\includegraphics[width=2.2in]{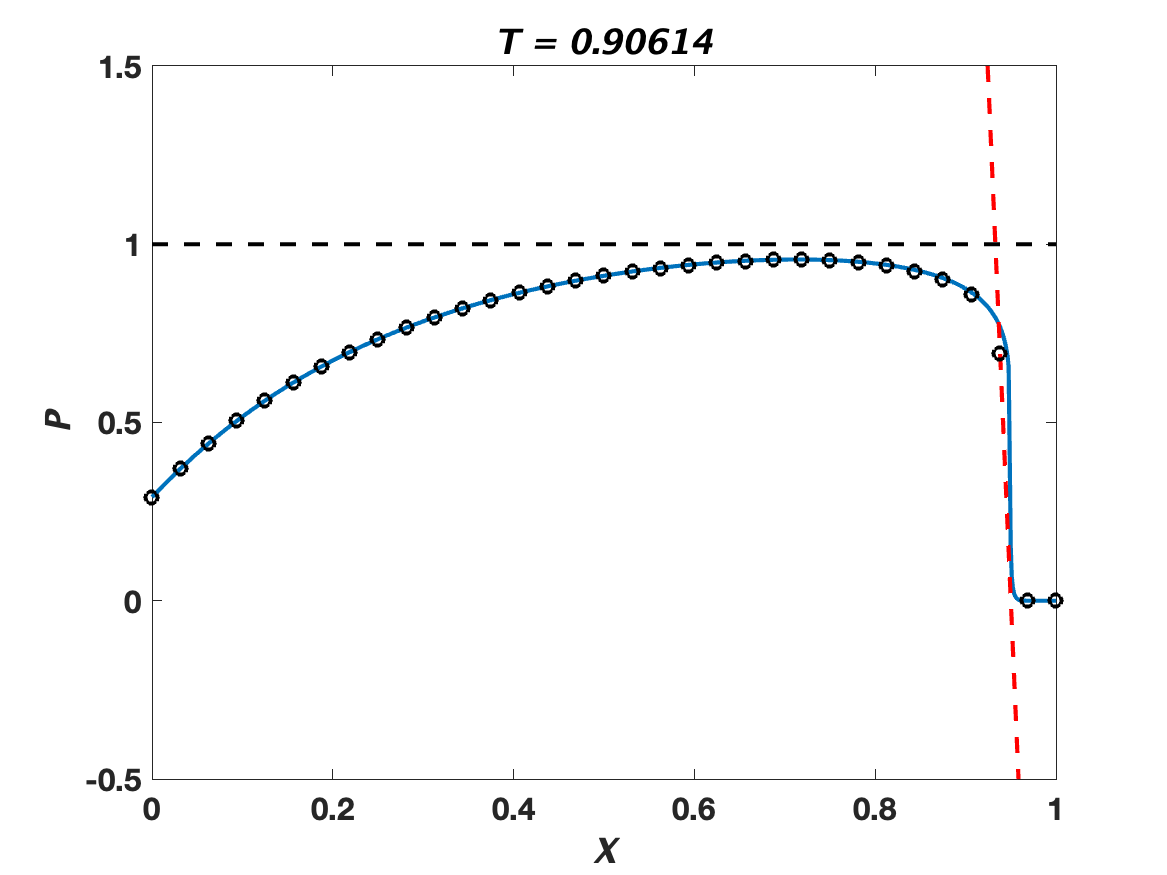} 
\caption{Blue curves:  Numerical solution of IBVP~\eqref{IBVP:LWE_short} based on the parameters given in 
Eq.~\eqref{eq:LWE1parms}, except here $\epsilon=0.4$, computed using the KSS method described in Section~\ref{sect:LWE_scheme} with  $N=262,144$ grid points and CFL number 10. The black circles represent 
the solution computed using FDS~\eqref{eq:FDS_P} with $M=4,096$  and CFL number $\bull \approx 0.4757$ (see Eq.~\eqref{eq:FDS_CFL}).  The red-dashed lines are plots of  $(X-\ram(T))\lshad P_{X} \rshad(T)$; they represent the tangent lines to  the solution profiles  at the wave-front $X=\ram(T)$.}
\label{figLWEm}
\end{figure}

Next, the sequence shown in Fig.~\ref{figLWEm} depicts the solution of the same
problem, with the same parameter values as in Eq.~\eqref{eq:LWE1parms}, except now
$\epsilon=0.4$; this increase in $\epsilon$, beyond the value of $\epsilon^{\bullet}$, yields
$\alpha=\alpha^{\bullet}\approx -0.200618$ and $T_{\rm f}=T_{\infty} = T_{1} \approx 0.951481$.  The solution is plotted at $T\approx 0.3, 0.6, 0.90614$, where $X = \ram (0.90614)\approx 0.95$.
The numerical solution is obtained using the same method, with the same discretization
and regularization parameters.
Again, it can be seen that the slopes of the KSS solution profiles at the wave-fronts match the analytically-derived values (see Eq.~\eqref{eq:[PX]}).  Furthermore, on comparing the former and latter figures we see  that $P$ exhibits growth (resp.~damping) for $\alpha >0$ (resp.~$\alpha <0$); this behavior is not unexpected, based on the results presented in Ref.~\cite[\S\,5.4]{PMJ22}, because the sign of the coefficient of $P_{X}$ in Eq.~\eqref{eq:EoM_P} is the same as $\sgn(\alpha)$.

\begin{figure}[ht]
\includegraphics[width=2.2in]{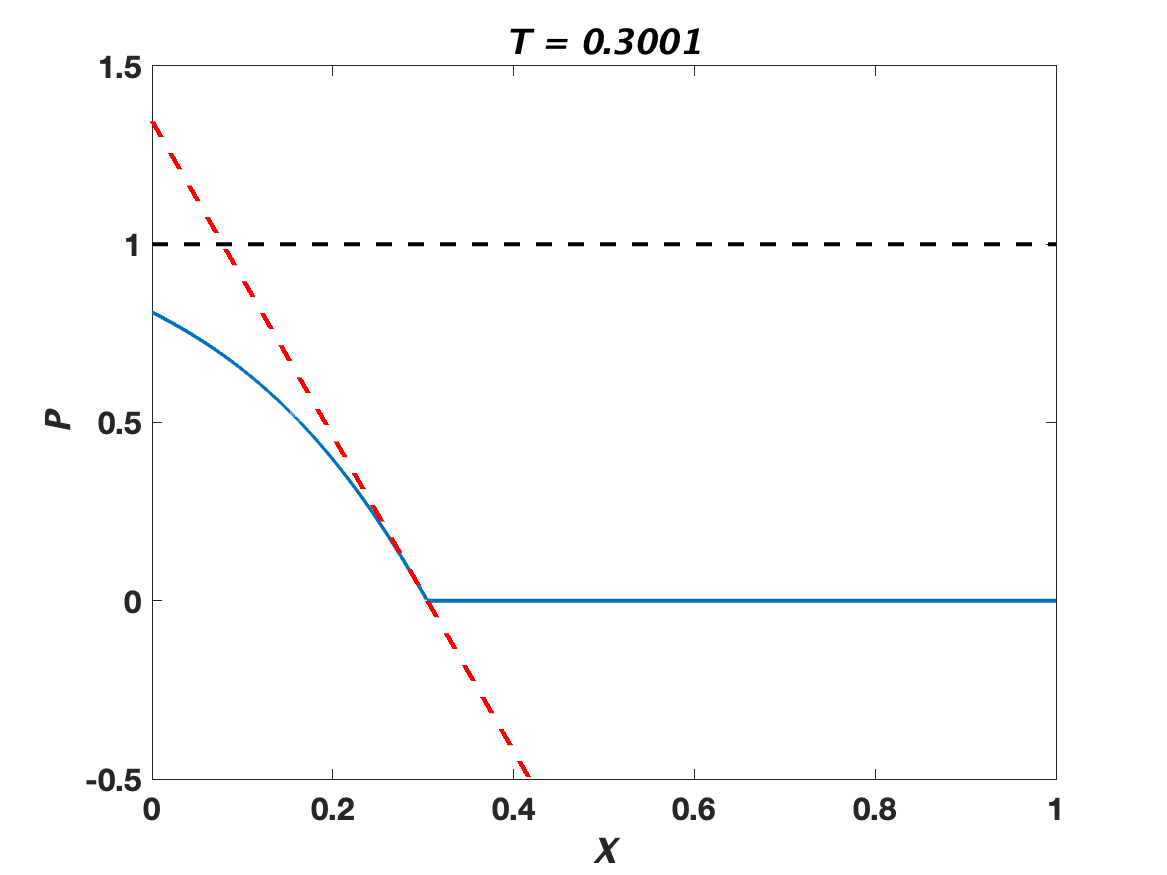}
\includegraphics[width=2.2in]{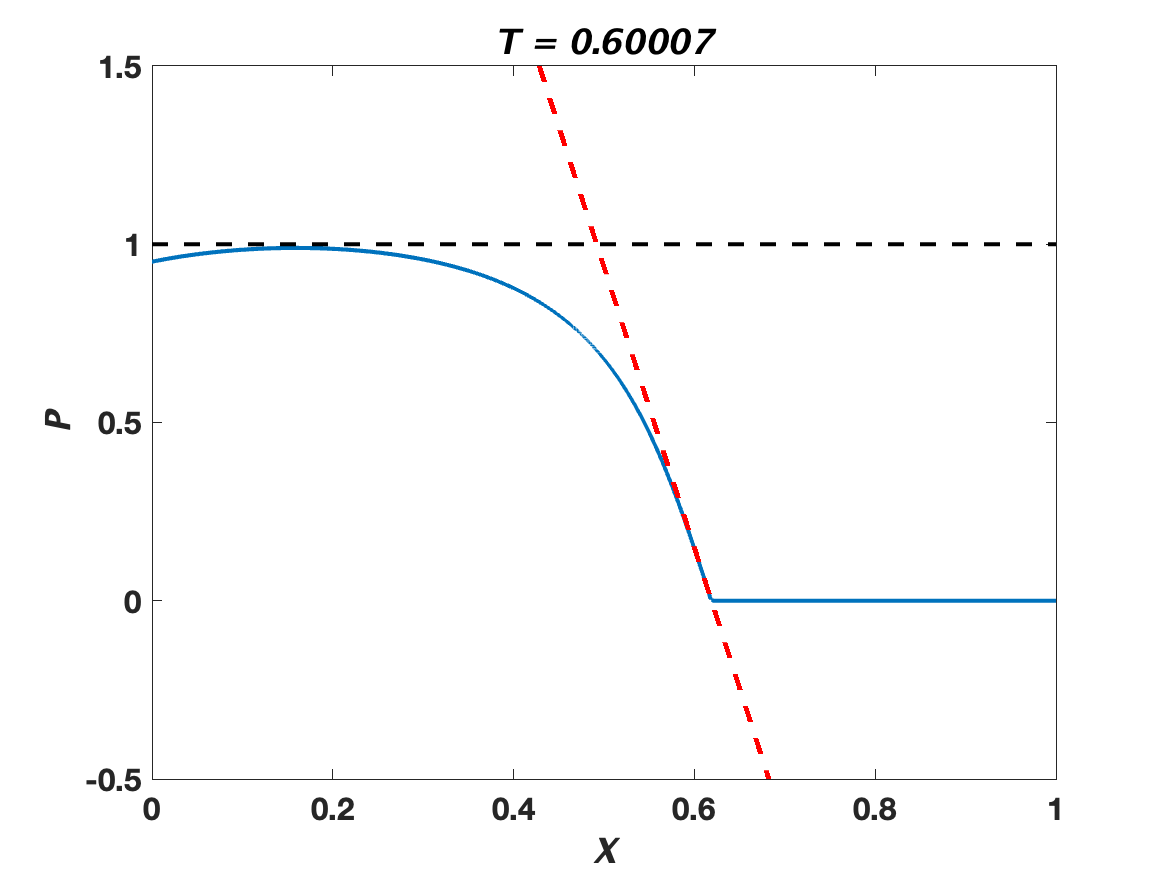}
\includegraphics[width=2.2in]{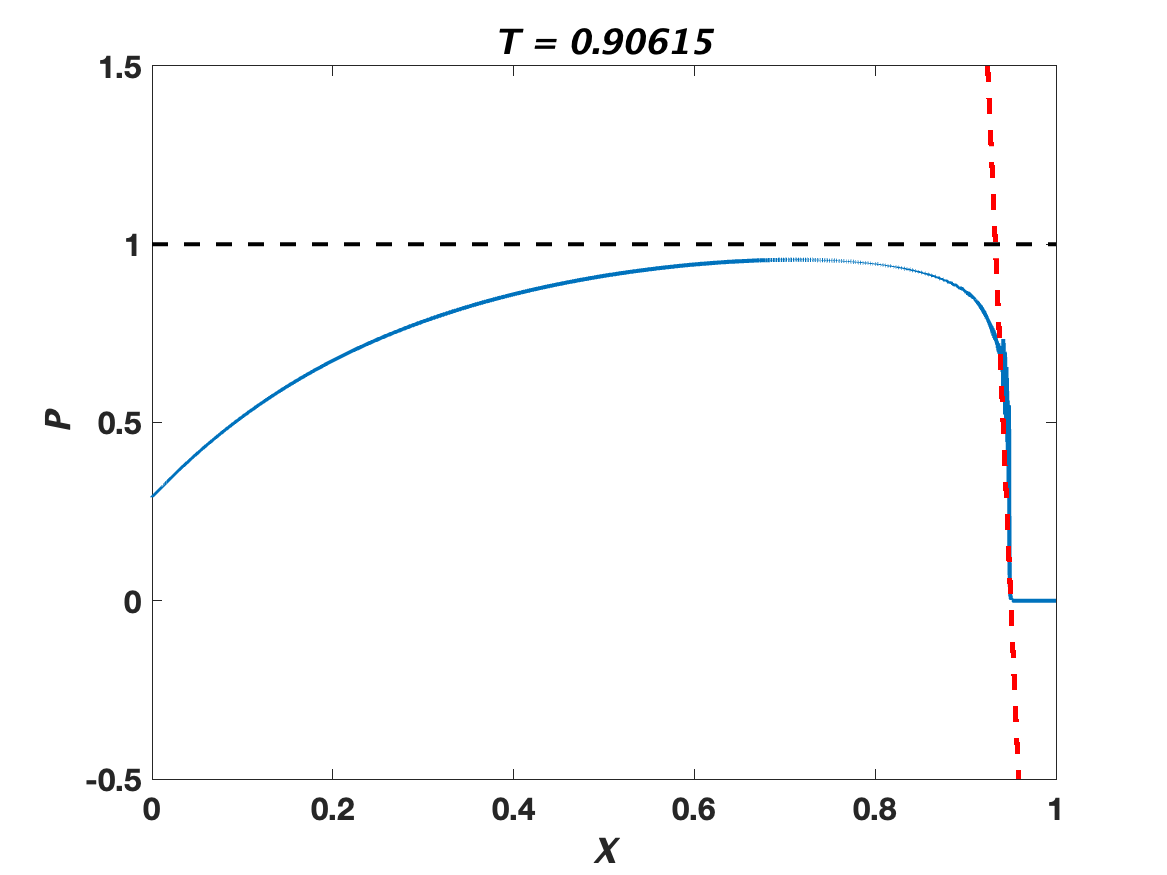}
\caption{Blue curves:  Numerical solution of  IBVP~\eqref{IBVP:LWE_short} using the same  
 parameters values as in Fig.~\ref{figLWEm}; theses solution curves were plotted from data sets generated by FDS~\eqref{eq:FDS_P} with $M =4,096$ and CFL number $\bull \approx 0.4757$ (see Eq.~\eqref{eq:FDS_CFL}). Again, the red-dashed lines are plots of  $(X-\ram(T))\lshad P_{X} \rshad(T)$; they represent the tangent lines to  the solution profiles  at the wave-front $X=\ram(T)$.}
\label{figLWEmFD}
\end{figure}

In concluding this subsection we call attention to the following:  While they show excellent agreement for the $T$-values taken in the left and center panels, in the right panel of Fig.~\ref{figLWEm} we see a clear discrepancy between the  KSS solution profile and the (selected) points from a data set generated by FDS~\eqref{eq:FDS_P} near $(\ram(T_{\rm end}),T_{\rm end})$.  The likely explanation for this lack of agreement can be seen  in the right panel of Fig.~\ref{figLWEmFD}, wherein the tangent line has been omitted for clarity.  There, one observes  that  the  FDS~\eqref{eq:FDS_P}-based solution profile, which is the result of plotting  all points in the aforementioned  data set,  exhibits high-frequency oscillations near $(\ram(T_{\rm end}),T_{\rm end})$, i.e., in exactly the same location where the disagreement occurs.

\subsection{Comparison with Krylov subspace methods}

To conclude this section, we modify our method for solving Eq.~\eqref{IBVP:LWE_short} by 
using Lanczos iteration to compute any required matrix function-vector products.  
Recall from Section \ref{subsubsec:LWEnumKSS}
that $\bar{\cal L}_N$ is the $N\times N$ matrix that is the spatial discretization of
the differential operator $\bar{\cal L}$ from Eq.~\eqref{eq:barLdef}.
To approximate $f(\bar{\cal L}_N){\bf u}$ for some $N$-vector ${\bf u}$ and function $f$, we
first apply $K$ iterations of
the Lanczos algorithm to $\bar{\cal L}_N$ with initial vector ${\bf u}$, which yields an orthonormal
basis $\{ {\bf q}_1, \ldots, {\bf q}_K \}$ for the $K$-dimensional Krylov subspace 
$${\cal K}({\bf u},\bar{\cal L}_N,K) = \textrm{span} \{ {\bf u}, \bar{\cal L}_N {\bf u},
\bar{\cal L}_N^2 {\bf u}, \ldots, \bar{\cal L}_N^{K-1} {\bf u} \},$$
and a symmetric tridiagonal matrix $T_K$ such that
$$\bar{\cal L}_N Q_K = Q_K T_K + \beta_K {\bf q}_{K+1} {\bf e}_K^T,$$
where $Q_K = \left[ \begin{array}{ccc} {\bf q}_1 & \cdots & {\bf q}_K \end{array} \right]$
and $\beta_K$ is a scalar.  Then, we have
\begin{equation} \label{eq:KSfunc}
f(\bar{\cal L}_N){\bf u} \approx {\bf w}_K = \|{\bf u}\|_2 Q_K f(T_K) {\bf e}_1.
\end{equation}
This approach was analyzed in Ref.~\cite{HL97} for the case of $f(\lambda)=e^{-\lambda\tau}$.
In our implementation, the number of iterations $K$ is chosen so that convergence of the
above approximation is achieved to within some tolerance; more precisely, the relative
difference $\|{\bf w}_{K+1}-{\bf w}_K\|_2/\|{\bf w}_{K+1}\|_2$ is required to be less than $10^{-4}$.

\begin{figure}[ht]
\includegraphics[width=2.2in]{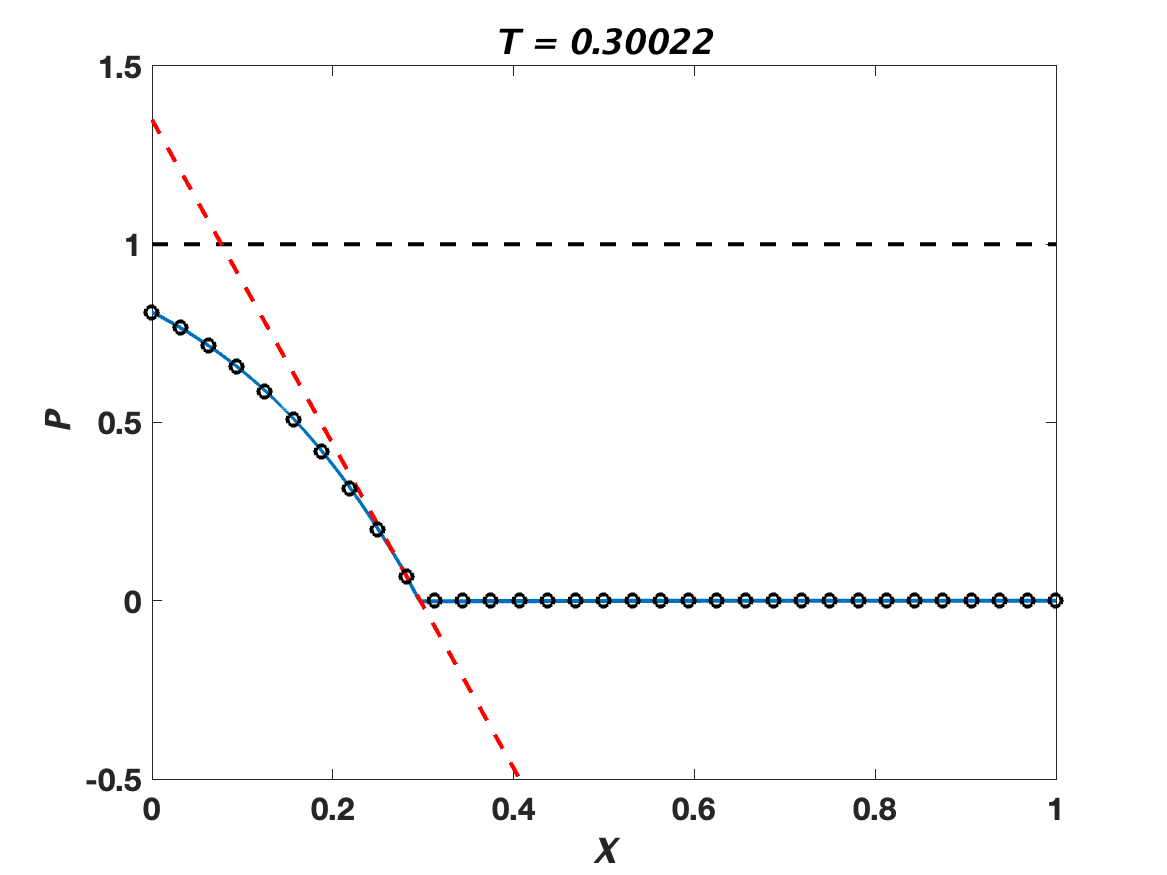}
\includegraphics[width=2.2in]{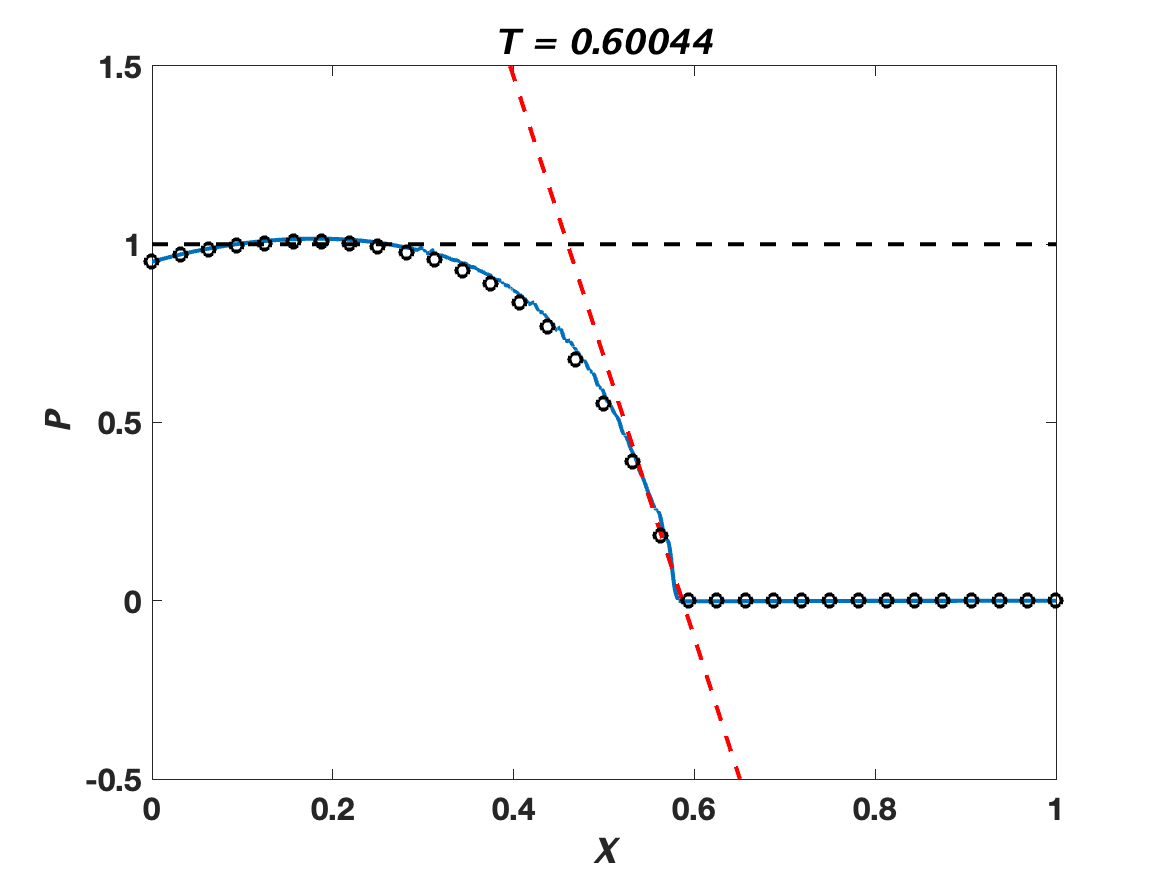}
\includegraphics[width=2.2in]{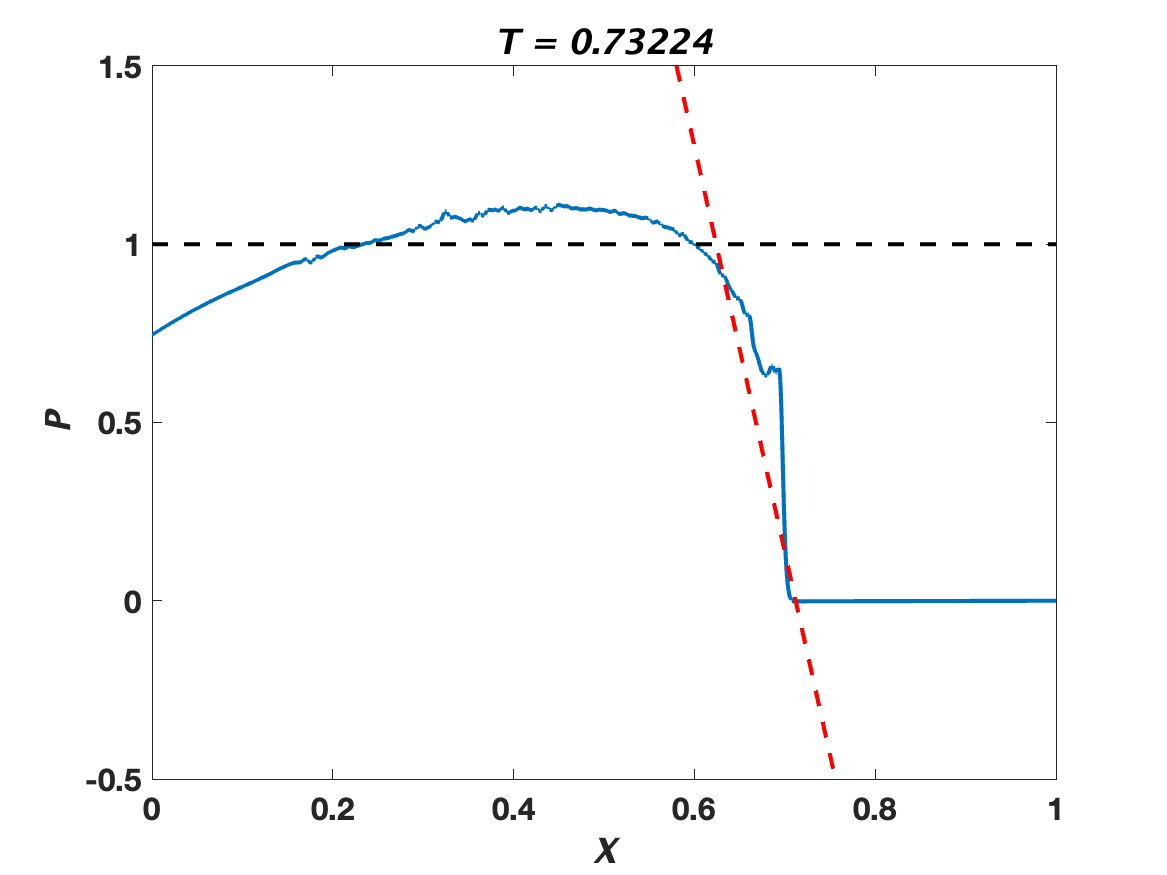}
\caption{Blue curves: Numerical solution of IBVP~\eqref{IBVP:LWE_short}  based on the parameters given in Eq.~\eqref{eq:LWE1parms}, computed using the numerical method described in Section~\ref{sect:LWE_scheme} except with matrix function-vector products approximated
as in Eq.~\eqref{eq:KSfunc}, using 
$N=4,096$ grid points and CFL number 1. The black circles represent the solution computed using FDS~\eqref{eq:FDS_P} with $M=4,096$  and CFL number $\bull \approx 0.5201$ (see Eq.~\eqref{eq:FDS_CFL}).  The red-dashed lines are plots of  $(X-\ram(T))\lshad P_{X} \rshad(T)$; they represent the tangent lines to  the solution profiles  at the wave-front $X=\ram(T)$.}
\label{figLWEKS}
\end{figure}

Results are presented in Fig.~\ref{figLWEKS}.  We use the same parameter values
as in Eq.~\eqref{eq:LWE1parms}.  The blue curves were obtained using $N=4,096$ grid
points and a CFL number of $\frac{\Delta T}{\Delta X}=10$.  Oscillations due to the Gibbs
phenomenon are regularized as they were when using the KSS method, except that
regularization is only applied to $P_{TT}$, and the Lanczos sigma factors are raised to
the power $\varpi^{\diamond\diamond}=4,096$.  We see that the solution is quite accurate
at $T=0.30022$, but loses accuracy and exhibits some small oscillations at $T=0.60044$.
The method breaks down at $T=0.73224$, as the Lanczos iteration fails to converge to within
the prescribed tolerance.
Without the aforementioned regularization, breakdown occurs even earlier.
Prior to breakdown, the Lanczos iteration required 11--15 matrix-vector products
for each matrix function-vector product, in contrast to the KSS method, which required only one matrix-vector
product and three Fourier transforms.

We also considered the application of an adaptive Krylov iteration from Ref.~\cite{adapkrylov}
to compute the matrix function-vector products in the spatial discretization of Eq.~\eqref{eq:LWEsolnoperator}.
The resulting time-stepping scheme is the exponential Euler method~\cite{ostermann22},
$$
\left[ \begin{array}{c} {\bf U}_N^{n+1} \\ {[{\bf U}_N']}^{n+1} \end{array} \right] = \phi_0(J_N\Delta T) 
\left[ \begin{array}{c} {\bf U}_N^{n} \\ {[{\bf U}_N']}^{n} \end{array} \right] + 
\Delta t\phi_1(J_N\Delta T) \left[ \begin{array}{c} {\bf 0} \\ {\bf b}_N^n \end{array} \right],
\quad J_N = \left[ \begin{array}{cc}
0 & I_N \\
-\bar{\cal L}_N & 0 
\end{array} \right],
$$
where $I_N$ is the $N\times N$ identity matrix and
${\bf b}_N^n$ is the spatial discretization of $b(Y,T_n)$ from Eq.~\eqref{eq:LWEbdef}.
Unfortunately, this approach fails to produce a reasonably accurate solution, even at a small 
time $T_n$, due to lack of smoothness of the solution.


\section{Closure}\label{sect:Close}


Extending work that began in Ref.~\cite{nmpde}, we have further extended the applicability
of KSS methods to wave propagation problems that include features such as Rayleigh dissipation,
and nonlinear effects that lead to finite-time shock formation. 
The result is a common framework within which KSS methods can accurately model a variety of wave phenomena, without having to satisfy a CFL condition of the classical type.  It has also been demonstrated that KSS methods are capable of producing accurate solutions of problems for which more established methods---Fourier spectral methods in the linear constant-coefficient case, and exponential integrators
using Krylov projection in the variable-coefficient or nonlinear case---break down.

Future work will include adding adaptivity to KSS methods, in which residual correction can be used
to not only automatically adjust the time step, as in Ref.~\cite{dozier}, but also the spatial resolution
and any regularization, to prevent contamination of the solution by high-frequency oscillations 
resulting from discontinuities while also avoiding excessive smoothing.  
The KSS framework presented in this paper will also be applied to other  inhomogeneous media propagation problems of interest, such as those involving dusty gases, poroacoustic drag, and media with multiple layers.

\section*{Acknowledgments}
B.R.\ was supported by NASA funding.  J.V.L.\ and P.M.J.\ were supported by 
U.S.~Office of Naval Research (ONR) funding. 

\appendix\markboth{}{}
\renewcommand{\thesection}{A}
\numberwithin{equation}{section}
\section{Appendix: Results from singular surface theory}\label{App:SST}

Adopting the usual notation convention, we let   $\lshad \mathfrak{F} \rshad (t)$ denote  the amplitude of the  jump discontinuity in the function  $\mathfrak{F}=\mathfrak{F}(\xi,t)$ across the singular surface $\xi=\Sigma(t)$; here,  we define $\lshad \mathfrak{F} \rshad (t)$ by
\begin{equation}\label{eq:jump-def}
\lshad \mathfrak{F} \rshad (t) :=\mathfrak{F}^{-}-\mathfrak{F}^{+},
\end{equation}
where $\mathfrak{F}^{\mp} := \lim_{\xi \to \Sigma (t)^{\mp}}\mathfrak{F}(\xi, t)$ are assumed to exist, and a `$+$' superscript corresponds to the region into which $\Sigma$ is advancing while a `$-$' superscript corresponds to the region behind $\Sigma$. 

In this communication,  $\xi=\Sigma (t)$ is a smooth planar surface  propagating along, and perpendicularly to,  the $\xi$-axis of our Cartesian coordinate system with velocity $V=V(t)$, where $V(t)=\rd \Sigma(t)/\rd t$.  If $\lshad \mathfrak{F} \rshad (t) \neq 0$ and $\mathfrak{F}$ represents a field variable (e.g., $\rho$), then $\xi=\Sigma (t)$ is termed a \emph{shock wave}~\cite[\S\,54]{Serrin59}; if, on the other hand, $\lshad \mathfrak{F} \rshad (t) =0$ but $\lshad \mathfrak{F}_{\xi} \rshad (t), \lshad \mathfrak{F}_{t} \rshad (t) \neq 0$, again assuming that  $\mathfrak{F}$ represents a field variable, then $\xi=\Sigma (t)$ is termed an \emph{acceleration wave}\footnote{For other studies wherein acoustic, or acoustic-like, versions of this class of singular surfaces arise, see, e.g., Refs.~\cite{BS22,JAMM67,JSV05,KJC18,Morro80,Sacc94,S08,S13,Walsh69}, and those cited therein.}.
 
Lastly, in Section~\ref{sect:ihLWE} we make use of the following tools from singular surface theory:  The 1D version of \emph{Hadamard's  lemma} (see, e.g., Refs.~\cite{Bland88}, \cite[p.~301]{S08}, and~\cite[\S\,174]{TT60}), specifically,
\be\label{eq:Hadamard}
\frac{\mathfrak{d} \lshad \mathfrak{F} \rshad (t)}{\mathfrak{d} t} = \lshad 
\mathfrak{F}_{t} \rshad(t) + V\lshad \mathfrak{F}_{\xi} \rshad(t), 
\en
where $\mathfrak{d}/\mathfrak{d} t$  gives the time-rate-of-change measured by an observer traveling with $\Sigma(t)$~\cite{Bland88};  the (1D) \emph{Maxwell compatibility condition}~\cite[\S\,175]{TT60}
\begin{equation}\label{eq:Max_CC}
\lshad 
\mathfrak{F}_{t} \rshad(t) + V(t)\lshad \mathfrak{F}_{\xi}\rshad(t) =0,
\end{equation}
which follows from Hadamard's lemma when $\lshad \mathfrak{F} \rshad (t)=0$; and the jump product formula~\cite[p.~302]{S08} 
\begin{equation}\label{eq:jump-prod}
\lshad \mathfrak{F}\mathfrak{G} \rshad= \mathfrak{F}^{+}\lshad \mathfrak{G} \rshad+\mathfrak{G}^{+}\lshad \mathfrak{F} \rshad+\lshad \mathfrak{F} \rshad \lshad \mathfrak{G} \rshad.
\end{equation}

\appendix\markboth{}{}
\renewcommand{\thesection}{B}
\numberwithin{equation}{section}
\section{Appendix: Simple explicit scheme for the ihLWE}\label{App:Scheme}

Generalizing the scheme developed in Ref.~\cite{JSV05}, wherein acoustic acceleration waves under the constant-coefficient (i.e., homogeneous fluid) version of the LWE were examined,  we construct the following simple discretization of Eq.~\eqref{eq:EoM_P}:
\begin{multline}\label{eq:LWE_Discrete}
\left(1- 2\epsilon \hat{\beta}  P_{m}^{k} \right) \!\left[\frac{P_{m}^{k+1}- 2P_{m}^{k}+P_{m}^{k-1}}{(\Delta T)^2} \right] 
-\, \exp(-\alpha X_{m})\!\left[\frac{P_{m+1}^{k}- 2P_{m}^{k}+P_{m-1}^{k}}{(\Delta X)^2}\right]\\ 
+\,\alpha \exp(-\alpha X_{m})\!\left[\frac{P_{m+1}^{k}- P_{m-1}^{k}}{2(\Delta X)}\right] 
= 2\epsilon \hat{\beta} \left(\frac{P_m^{k}- P_m^{k-1}}{\Delta T}\right)^{2},
\end{multline}
where we recall that
\be
\hat{\beta} := \beta/\gamma = (1+1/\gamma/2).
\en
In Eq.~\eqref{eq:LWE_Discrete}, $P(X_m,T_k)\approx P_m^k$, where the mesh points $(X_m, T_k)$ are given by $X_m=m(\Delta X)$ and $T_k=k (\Delta T)$. Additionally, and as dictated by  IBVP~\eqref{IBVP:LWE_non-D}, the spatial- and temporal-step sizes are defined as $\Delta X=1/M$ and $\Delta T=T_{\rm f}/(2M)$, respectively, where $M\gg 1$ is an integer.
 
On solving  for  $P_{m}^{k+1}$, the most advanced time-step approximation, we obtain the explicit finite difference scheme~(FDS)
\begin{multline}\label{eq:FDS_P}
P_{m}^{k+1} =2P_{m}^{k}-P_{m}^{k-1}+ \left(1- 2\epsilon \hat{\beta} P_{m}^{k} \right)^{-1}
\Bigg{[}\bull^{2}\exp(-\alpha X_{m})\!\left(P_{m+1}^{k}- 2P_{m}^{k}+P_{m-1}^{k} \right)\\ 
-\,\alpha \bull (\Delta T) \exp(-\alpha X_{m})\!\left(P_{m+1}^{k}- P_{m-1}^{k}\right)/2
+ 2\epsilon \hat{\beta}  \left(P_m^{k}- P_m^{k-1}\right)^{2}\Bigg{]},
\end{multline}
which holds for each $m=1,2,3, \ldots, M-1$ and each $k=1,2,3, \ldots, 2M-1$.  Here, we have set
\be \label{eq:FDS_CFL}
\bull := \frac{\Delta T}{\Delta X} =\tfrac{1}{2}T_{\rm f} < 1,
\en
which of course is also the CFL number of our (dimensionless) FDS.
 
And lastly, the discretized versions of the  BCs and ICs of IBVP~\eqref{IBVP:LWE_non-D} read
\be
P_{0}^{k} =\sin(\pi T_{k}), \qquad \qquad P_{M}^{k}=0 \qquad (k=0,2,3, \ldots, 2M)
\en
and
\be
P_{m}^{0} = 0, \qquad \qquad P_{m}^{1}= P_{m}^{0} \qquad (m=1,2,3, \ldots, M-1),
\en
respectively. 


\end{document}